\documentclass[article,11pt]{aiaa-pretty}    

\usepackage{amssymb}
\usepackage{amsmath}
\usepackage{booktabs}
 \usepackage{varioref}
 \usepackage{wrapfig}
 \usepackage{threeparttable}
 \usepackage{dcolumn}
  \newcolumntype{d}{D{.}{.}{-1}}
 \usepackage{nomencl}
  \makeglossary
 \usepackage{subfigure}
 \usepackage{subfigmat}
 \usepackage{fancyvrb}
  \fvset{fontsize=\footnotesize,xleftmargin=2em}
 \usepackage{lettrine}

\newcommand{\Real}{\mathbb R}

\newcommand{\real}[1]{{\mathbb R}^{#1}}
\newcommand{\Ebar}{{\overline{E}}}

\newcommand{\be}{{\boldsymbol e}}
\newcommand{\bff}{{\boldsymbol f}}

\newcommand{\bh}{{\boldsymbol h}}

\newcommand{\bq}{{\boldsymbol q}}
\newcommand{\br}{{\boldsymbol r}}

\newcommand{\bu}{{\boldsymbol u}}

\newcommand{\bx}{\boldsymbol x}

\newcommand{\bxf}{{\bx(\cdot)}}  
\newcommand{\bxt}{{\widetilde{\bx}}}
\newcommand{\buf}{{\bu(\cdot)}}  
\newcommand{\but}{{\widetilde{\bu}}}  

\newcommand{\bP}{{\boldsymbol P}}

\newcommand{\X}{\mathbb{X}}
\newcommand{\U}{\mathbb{U}}

\newcommand{\bzero}{{\bf 0}}

\newcommand{\bnu}{{\mbox{\boldmath $\nu$}}}
\newcommand{\bmu}{{\mbox{\boldmath $\mu$}}}

\newcommand{\blam}{{\mbox{\boldmath $\lambda$}}}

\newcommand{\tlam}{{\mbox{\boldmath $\widetilde{\lambda}$}}}
\newcommand{\tmu}{{\mbox{\boldmath $\widetilde{\mu}$}}}
\newcommand{\tnu}{{\mbox{\boldmath $\widetilde{\nu}$}}}

\newcommand{\bet}{{\widetilde{\be}}}
\newcommand{\bht}{{\widetilde{\bh}}}
\newcommand{\blamt}{{\widetilde{\blam}}}
\newcommand{\bmut}{{\widetilde{\bmu}}}
\newcommand{\bnut}{{\widetilde{\bnu}}}
\author{ %
I. M. Ross\thanks{Distinguished Professor and Program Director, Control and Optimization, Department of Mechanical and Aerospace Engineering.}\\
\textit{Naval Postgraduate School, Monterey, CA 93943}\\[1em]
Q. Gong\thanks{Associate Professor, Department of Applied Mathematics.}\\
\textit{University of California, Santa Cruz, CA 95064}\\[1em]
M. Karpenko\thanks{Research Associate Professor and Associate Director, Control and Optimization Laboratories, Department of Mechanical and Aerospace Engineering; corresponding author: \texttt{mkarpenk@nps.edu}.},
R. J. Proulx\thanks{Research Professor, Control and Optimization Laboratories, Space Systems Academic Group.}\\
\textit{Naval Postgraduate School, Monterey, CA 93943}
}
\title{Scaling and Balancing for High-Performance Computation of Optimal Controls}

\abstract{
It is well-known that proper scaling can increase the efficiency of  computational problems. In this paper we define and show that a balancing technique can substantially improve the computational efficiency of optimal control algorithms.  We also show that non-canonical scaling and balancing procedures may be used quite effectively to reduce the computational difficulty of some hard problems. These results have been used successfully for several flight and field operations at NASA and DoD.  A surprising aspect of our analysis shows that it may be inadvisable to use auto-scaling procedures employed in some software packages.  The new results are agnostic to the specifics of the computational method; hence, they can be used to enhance the utility of any existing algorithm or software.
}

\begin{document}
\maketitle


\section{Introduction}

In many practical optimal control problems, the decision variables range ``wildly'' in several orders of magnitude\cite{conway:survey}. For instance, in a space trajectory optimization problem\cite{conway:book-chapter,longuski}, a position variable can vary from a few meters to well over a million kilometers.   The conventional wisdom to manage the associated computational problems is to use canonical units.  Many space trajectory optimizations fare well when canonical units are used\cite{conway:book-chapter, longuski}. For instance, for low-Earth-orbiting spacecraft, one may choose the radius of the Earth ($R_\oplus$) as a unit of distance, and circular speed (at $R_\oplus$) as the unit of velocity.  For interplanetary spacecraft, astronomical units provide a set of canonical units for scaling trajectory optimization problems.
\textbf{\emph{In this paper, we show that
it is possible to define and use arbitrary and inconsistent units for faster trajectory optimization}}.  For example, we may choose meters as a unit of distance along the $x$-direction while concurrently using feet for distance units along the $y$-direction. Furthermore, we show that it is not necessary to choose a consistent (or canonical) unit of velocity in the $x$-direction to be equal to meters per second (or feet per second for velocity in the $y$-direction).  Thus, for example, one may ``arbitrarily'' choose yards per day as the unit of $x$-velocity while insisting that the $x$-position be measured in meters. The purpose of using such unusual or \textbf{\textrm{\emph{designer units}}}\cite{ross-book}; i.e., highly customized units that do not necessarily conform to standardized units,  is to liberate ourselves from using well-established canonical/consistent units so that we may scale an optimal control problem for faster computational results.  This liberation allows us to radically alter what we mean by scaling optimal control problems, and consequently solve some apparently ``hard'' problems with more ease than ever before.

A second aspect of our paper is the impact of scaling on dual variables. From the Hahn-Banach theorem\cite{Tao:epsilon}, dual variables exist for all computational optimal control problems even when they are not used. To illustrate this consequential theorem, consider the ordinary differential equation,
\begin{equation} \label{eq:ode}
\dot\bx = \bff(\bx)
\end{equation}
where, $ \bx \in \real{N_x}$ and  $\bff: \bx \mapsto \real{N_x}$.
The formal adjoint to the variation of \eqref{eq:ode} is defined by
\begin{equation} \label{eq:adjoint}
\dot\blam = -\left[\partial_x\bff\right]^T\, \blam
\end{equation}
where, $\partial_x\bff$ is the Jacobian of $\bff$ with respect to $\bx$.
\emph{\textbf{In other words, \eqref{eq:adjoint} exists from the mere fact that \eqref{eq:ode} exists}}; hence, when a differential equation in a computational optimal control problem is scaled, it automatically affects the adjoint equation, even if it (i.e. the adjoint equation) is not used.  In this paper, we show that the equations in a computational optimal control problem must be scaled in such a way that it does not ``unscale'' the adjoint variable even if the adjoint equation is never used in the algorithm. We call this type of scaling ``balancing.''

In most algorithms -- including the ones where derivatives are not used (e.g., genetic algorithms) --  the ``information content'' in the Jacobian forms a key ingredient in the recipe that connects the sequence of iterations\cite{McGrath:napa-genetic,McGrath:napa-unscent,McGrath:diss}. Consequently, the scales used in an optimal control problem must be balanced even if the adjoint equation is never used because it represents the information content contained in a Jacobian by way of~\eqref{eq:adjoint}.
In other words, there is no escape from considering the adjoint equations.  This is a fundamental result traceable to the Hahn-Banach theorem.  An alternative explanation for the no-escape clause is that the adjoint equations are part of the necessary conditions for optimality.  By definition, necessary conditions are indeed necessary; hence, it should not be entirely surprising that balancing is also necessary.

Our analysis also reveals another surprising result: that if scaling is performed at the discrete level, it inadvertently introduces new terms in the dynamical equations with possible feedback effects that may destabilize the search algorithm.  Consequently, \emph{\textbf{automatic techniques that scale the problem at the discrete level may be more harmful than useful}}. The simple remedy is to scale and balance the equations at the optimal-control level and choose algorithms that do not scale the equations at the discrete level.  This simple ``trick'' has been used many number of times before by NASA\cite{bhatt:opm,SIAMnews,zpm:NASA-report,zpm:IEEE,TEI-JGCD-2011,TRACE-IEEE-Spectrum,Kepler-micro-slew,Karp-JWST} and DoD\cite{CDC-Workshop,TRACE-JGCD-2014,RobStevens-Tether,PSReview-ARC-2012,Minelli-AAS,RS-Expt} to solve and implement on-orbit and fielded solutions.  Note, however, that scaling and balancing were not at the forefront in many of these applications because their main contributions far outweighed such discussions.  For instance, the focal point of \cite{bhatt:opm} was the flight implementation of the optimal propellant maneuver onboard the International Space Station rather than the employment of scaling and balancing techniques.  Similarly, \cite{zpm:NASA-report} and \cite{zpm:IEEE} were focused on the first flight implementations of a historic zero-propellant maneuver while \cite{Kepler-micro-slew} was on the feasibility of arcsecond slews for precision pointing of the Kepler spacecraft.  In the same spirit, the main contribution of \cite{TRACE-JGCD-2014} was the flight implementation of a shortest-time maneuver, and not on the specifics or the importance of scaling and balancing that were necessary to accomplish the on-orbit demonstration.

In this paper, we generalize the application-specific procedures of the past successes by laying down the mathematical foundations for scaling and balancing of generic optimal control algorithms.  In doing so,\textbf{\emph{ we also demonstrate the fallacies of some scaling techniques that are currently in practice}}.

\vskip 0.4in
\fbox{
\begin{minipage}{0.9\textwidth}
The journal version of this paper contains some typographical errors.  Regardless, please cite the journal paper (\textit{J. Guid., Contr. \& Dyn., 41/10, 2018, pp.~2086--2097}) if you need to cite the results contained herein.
\end{minipage}
}

\section{General Problem Formulation}
A generic optimal control problem can be formulated as follows:
%
\newlength{\leftside}
\newlength{\rightside}
\newcommand*{\leftterm}{}
\newcommand*{\rightterm}{}
\newcommand*{\term}[1]{$\displaystyle#1$}
\begin{align*}
&\left.\begin{aligned}
\phantom{preamble}
\X & = \real{N_x}   &\U & = \real{N_u}& \\
\bx &= (x_1, \ldots, x_{N_x})  \quad &\bu &= (u_1, \ldots, u_{N_u}) & \\
\end{aligned}\hspace{2.85cm} \right\} \ \text{\emph{\textbf{preamble}}}
\\
%
&\begin{aligned}
\renewcommand*{\leftterm}{\phantom{\quad J[\bx(\cdot), \bu(\cdot), t_0, t_f]}}
\renewcommand*{\rightterm}{\phantom{\be(\bx_0, \bx_f, t_0, t_f) \le \be^U dt }}
\settowidth{\leftside}{\term{\leftterm}}
\settowidth{\rightside}{\term{\rightterm}}
\overbrace{(B)}^{\text{\normalsize \emph{\textbf{problem}}}} \left\{
\begin{array}{ll}
\text{Minimize }& \\[-1.25em] 
& \left.\begin{aligned}
\makebox[\leftside][r]{\term{J[\bx(\cdot), \bu(\cdot), t_0, t_f] := }} & \\
 \makebox[\leftside][r]{\term{E(\bx_0, \bx_f, t_0, t_f)}} & + \makebox[\rightside][l]{\term{\int_{t_0}^{t_f}F(\bx(t), \bu(t), t)\ dt }}  \\
 \end{aligned}\right\} \ \text{\textbf{\emph{cost}}}\\ [2em]
\text{Subject to} &  
\\[-1.5em]
&  \left.\begin{aligned}
    \makebox[\leftside][r]{\term{\dot\bx}} &= \makebox[\rightside][l]{\term{\bff(\bx(t), \bu(t), t)}}\\
  \end{aligned}\right\}\ \text{\emph{\textbf{dynamics}}} \\ [1em]
  %
&  \left.\begin{aligned}
    \makebox[\leftside][r]{\term{\be^L}}&\le \makebox[\rightside][l]{\term{\be(\bx_0, \bx_f, t_0, t_f) \le \be^U}}\\
  \end{aligned}\right\} \ \text{\emph{\textbf{events}}} \\[1em]
  &  \left.\begin{aligned}
   \makebox[\leftside][r]{\term{\bh^L}}
   &\le \makebox[\rightside][l]{\term{\bh(\bx(t), \bu(t), t) \le \bh^U}}\\
  \end{aligned}\right\} \ \text{\emph{\textbf{path}}}
\end{array}
\right.
\end{aligned}\label{prob:B}
\end{align*}
%

Although the symbols used in Problem $B$ are fairly standard (see, for example \cite{ross-book}), we note the following for the purposes of completeness:
\begin{itemize}
\item $\X$  and $\U$ are $N_x$- and $N_u$-dimensional real-valued state and control spaces respectively. We assume $N_x \in \mathbb{N}^+ $ and $N_u \in \mathbb{N}^+ $.
\item $J$ is the scalar cost function.  The arguments of $J$ are the optimization variables.
\item The optimization variables are:
\begin{itemize}
\item $\bxf$: the $N_x$-dimensional state trajectory,
\item $\buf$: the $N_u$-dimensional control trajectory,
\item $t_0$: the initial clock time, and
\item $t_f$: the final clock time.
\end{itemize}
\item $E$ is the scalar \emph{\textbf{endpoint cost function}}.  The arguments of $E$ are the endpoints.  In the classical literature, $E$ is known as the ``Mayer'' cost function.
\item The endpoints are the initial state $\bx_0 \equiv \bx(t_0)$, the final state $\bx_f \equiv \bx(t_f)$, the initial time $t_0$ and the final time $t_f$.
\item $F$ is the scalar \emph{\textbf{running cost function}}.  The arguments of $F$ are the instantaneous value of the state variable $\bx(t)$, the instantaneous value of the control variable $\bu(t)$ and time $t$.
\item $\bff$ is the $N_x$-dimensional ``\textbf{\emph{dynamics function}},'' or more appropriately the right-hand-side of the dynamics equation. The arguments of $\bff$ are exactly the same as the arguments of $F$.
\item $\be$ is the $N_e$-dimensional \emph{\textbf{endpoint constraint function}}.  The arguments of $\be$ are exactly the same as that of $E$.
\item $\be^L$ and $\be^U$ are the $N_e$-dimensional lower and upper bounds on the values of $\be$.
\item $\bh$ is the $N_h$-dimensional \textbf{\emph{path constraint function}}. The arguments of $\bh$ are exactly the same as that of $F$.
\item $\bh^L$ and $\bh^U$ are the $N_h$-dimensional lower and upper bounds on the values of $\bh$.
\end{itemize}
%
%
The five functions, $E$, $F$, $\bff$, $\be$ and $\bh$ are collectively known as the \emph{\textbf{data functions}} (for Problem $B$).

Regardless of any type of method used to solve Problem~$B$ -- including the so-called direct methods -- a solution to the problem must at least satisfy the necessary conditions of optimality because they are indeed necessary.  The necessary conditions for Problem $B$ are generated by an application of  Pontryagin's Principle\cite{ross-book,vinter,clarke-2013book}.  This results in the following boundary value problem (BVP), which we denote as Problem $B^\lambda$:
%
%
\begin{eqnarray*}\label{prob:Plam}
& (\textsf{$B^\lambda$}) \left\{
\begin{array}{lrll}
& \dot{\bx}(t)-\partial_\lambda \overline{H}(\bmu(t), \blam(t), \bx(t),\bu(t), t) = &\bzero &\text{\textbf{\emph{(state eqns)}}} \\[.5em]
& \dot{\blam}(t)+\partial_x \overline{H}(\bmu(t), \blam(t), \bx(t),\bu(t), t) = &\bzero &\text{\emph{\textbf{(costate eqns)}}} \\[.5em]
&\bh^L \le \bh(\bx(t), \bu(t), t) \le & \bh^U &\text{\emph{\textbf{(path constraint)}}}\\[.75em]
& \partial_u \overline{H}(\bmu(t), \blam(t), \bx(t),\bu(t), t) = &\bzero &\text{\textbf{\emph{(Hamiltonian}}}\\
&\bmu & \!\!\!\dagger\  \bh &\text{\quad \emph{\textbf{Minimization)}}}\\[.75em]
&\be^L \le \be\big(\bx_0, \bx_f, t_0, t_f \big) \leq &\be^U &\text{\emph{\textbf{(endpoint eqns)}}}\\[.75em]
&\blam(t_0) + \partial_{x_0}\overline{E}(\bnu,\bx_0,\bx_f, t_0, t_f) = &\bzero &\text{\emph{\textbf{(initial and}} }\\
&\blam(t_f) - \partial_{x_f}\overline{E}(\bnu,\bx_0,\bx_f, t_0, t_f) = &\bzero &\text{\emph{\textbf{final transversality}}}\\
&\bnu &\!\!\! \dagger\ \be &\text{\textbf{\emph{conditions)}}}\\[.75em]
& \mathcal{H}[@ t_0] - \partial_{t_0}\overline{E}(\bnu,\bx_0,\bx_f, t_0, t_f) = &0 &\text{\textbf{\emph{(Hamiltonian}}}\\
& \mathcal{H}[@ t_f] + \partial_{t_f}\overline{E}(\bnu,\bx_0,\bx_f, t_0, t_f) = &0 &\text{\emph{\textbf{value conditons)}}}\\
\end{array} \right. &
\end{eqnarray*}
%
%
%
In Problem $B^\lambda$, the unknowns are:
\begin{enumerate}
\item The system trajectory, $t \mapsto (\bx, \bu) \in \real{N_x} \times \real{N_u}$;
\item The adjoint covector function, $ t \mapsto \blam \in \real{N_x}$;
\item The path covector function, $ t \mapsto \bmu \in \real{N_h}$;
\item The endpoint covector, $\bnu \in \real{N_e}$; and
\item The initial and final clock times, $t_0 \in \Real$ and $t_f \in \Real$.
\end{enumerate}
The quantity $\overline{H}$ is the \textbf{\emph{Lagrangian of the Hamiltonian}} given by
$$\overline{H}(\bmu, \blam, \bx, \bu, t) :=  H(\blam, \bx, \bu, t) + \bmu^T\bh(\bx, \bu, t)  $$
where $H$ is the usual Pontryagin Hamiltonian,
$$ H(\blam, \bx, \bu, t) := F(\bx,
\bu, t) + \blam^T \bff(\bx,
\bu, t) $$
and $\mathcal{H}$ is the \textbf{\emph{lower or minimized Hamiltonian}},
$$\mathcal{H}(\blam, \bx, t) := \min_{\bu} H(\blam, \bx, \bu, t)  $$
The symbol $\mathcal{H}[@t]$ is a shorthand for $\mathcal{H}(\blam(t), \bx(t), t)$.
The quantity $\Ebar$ is the \emph{\textbf{Endpoint Lagrangian}} given by,
$$
\Ebar(\bnu, \bx_0, \bx_f, t_0, t_f) := E(\bx_0, \bx_f, t_0, t_f) + \bnu^T \be(\bx_0, \bx_f, t_0, t_f)$$
The $\dagger$ notation used in defining Problem $B^\lambda$  denotes \textbf{\emph{complementarity conditions}}\cite{ross-book}.  For instance, $\bnu\, \dagger\, \be$ is a shorthand for the conditions
\begin{equation}\label{eq:nu-perp}
\nu_i \left\{\begin{array}{ccrc}
               \le 0            & if & e_i(\bx_0, \bx_f, t_0, t_f) &= e_i^L \\
               =  0             & if & \qquad e_i^L < e_i(\bx_0, \bx_f, t_0, t_f) &< e_i^U \\
               \ge 0            & if & e_i(\bx_0, \bx_f, t_0, t_f) &= e_i^U \\
               unrestricted     & if & e_i^L &= e_i^U
             \end{array}
   \right.
\end{equation}
with $\bmu\, \dagger \, \bh$ defined similarly (for each $t$).

\section{Scaling the Primal Problem}

We change the coordinates of the unknown variables of Problem $B$ according to the affine transformations,
\begin{subequations}\label{eq:transform-var}
\begin{align}
\bx := & \bP_x\ \bxt + \bq_x\\
\bu := & \bP_u\ \but + \bq_u \\
t := & p_t\ \widetilde{t} + q_t
\end{align}
\end{subequations}
where, the uppercase letter $\bP_{(\cdot)}$ denotes an invertible square matrix of appropriate dimensions, and the lower case letter $p_{(\cdot)}$ is a scalar.  Similarly $\bq_{(\cdot)}$ is a vector of appropriate dimension and $q_{(\cdot)}$ is a scalar.  The tilde ($ \sim $) variables are the transformed variables.  In similar fashion, we ``scale'' the cost functional $J$, the endpoint constraint function $\be$ and the path constraint function $\bh$ according to,
\begin{subequations}\label{eq:transform-fun}
\begin{align}
J := & p_J\ \widetilde{J} + q_J\\
\be := & \bP_e\ \bet + \bq_e \\
\bh := & \bP_h\ \bht + \bq_h
\end{align}
\end{subequations}
where, $p_J > 0$, and $\bP_e$ and $\bP_h$ are positive definite diagonal matrices.
Let Problem $\widetilde{B}$ denote the transformation of Problem $B$ resulting from \eqref{eq:transform-var} and \eqref{eq:transform-fun}.  This problem can be explicitly obtained as follows: First, the transformation of $J$ to $\widetilde{J}$ can be constructed using \eqref{eq:transform-var} as,
\begin{align}
J[\bxf, \buf, t_0, t_f] &= J[\bP_x \bxt(\cdot) + \bq_x,\ \bP_u \but(\cdot) + \bq_u,\ p_t \widetilde{t}_0 +q_t,\ p_t \widetilde{t}_f + q_t] \nonumber\\
&:= p_J \widetilde{J}[\bxt(\cdot), \but(\cdot), \widetilde{t}_0, \widetilde{t}_f ] + q_J \label{eq:J-tilde-basic}
\end{align}
where the last equality in \eqref{eq:J-tilde-basic} follows from \eqref{eq:transform-fun}. Hence, the cost functional transforms according to,
\begin{align}
\widetilde{J}[\bxt(\cdot), \but(\cdot), \widetilde{t}_0, \widetilde{t}_f ] :=& - \frac{q_J}{p_J} + \frac{1}{p_J} J[\bP_x \bxt(\cdot) + \bq_x,\ \bP_u \but(\cdot) + \bq_u,\ p_t \widetilde{t}_0 +q_t,\ p_t \widetilde{t}_f + q_t] \nonumber \\
=& - \frac{q_J}{p_J} + \frac{1}{p_J} E(\bP_x \bxt_0 + \bq_x,\ \bP_x \bxt_f + \bq_x,\ p_t \widetilde{t}_0 +q_t,\ p_t \widetilde{t}_f + q_t) \nonumber \\
& + \frac{p_t}{p_J} \int_{\widetilde{t}_0}^{\widetilde{t}_f} F(\bP_x \bxt(t) + \bq_x,\ \bP_u \but(t) + \bq_u,\ p_t \widetilde{t} +q_t)\ d\widetilde{t} \label{eq:J-tilde}
\end{align}
where $t$ in the integrand in \eqref{eq:J-tilde} is to be understood as $(p_t \widetilde{t} +q_t)$.

The transformation of the dynamics is given by,
\begin{align*}
\frac{d\bx}{dt} = \frac{\bP_x}{p_t} \frac{d\bxt}{d\widetilde{t}} = \bff(\bP_x \bxt(t) + \bq_x,\ \bP_u \but(t) + \bq_u,\ p_t \widetilde{t} +q_t)
\end{align*}
Hence, the dynamical equations transform according to,
\begin{align}\label{eq:dynamics-tilde}
\frac{d\bxt}{d\widetilde{t}} = \widetilde{\bff}(\bxt(t), \but(t), \widetilde{t}) := p_t \bP^{-1}_x \bff(\bP_x \bxt(t) + \bq_x,\ \bP_u \but(t) + \bq_u,\ p_t \widetilde{t} +q_t)
\end{align}
where, we have once again used $t$ in \eqref{eq:dynamics-tilde} to mean $(p_t \widetilde{t} +q_t)$.

By using similar procedures, it follows that the endpoint and path constraint functions transform according to,
\begin{align}
\bet\big(\bxt_0, \bxt_f, \widetilde{t}_0, \widetilde{t}_f \big) = &
\bP_e^{-1}\left[\be\big(\bP_x \bxt_0 + \bq_x,\ \bP_x \bxt_f + \bq_x,\ p_t \widetilde{t}_0 +q_t,\ p_t \widetilde{t}_f + q_t  \big) - \bq_e\right] \label{eq:events-tilde}\\
\bht(\bxt(\widetilde{t}), \but(\widetilde{t}), \widetilde{t}) = & \bP_h^{-1}\big[\bh \big(\bP_x \bxt(t) + \bq_x,\ \bP_u \but(t) + \bq_u,\ p_t \widetilde{t} +q_t\big) - \bq_h \big] \label{eq:path-tilde}
\end{align}
The corresponding lower and upper bounds are given by,
\begin{align}
\widetilde{\be}^L &=  \bP_e^{-1}(\be^L - \bq_e), &\widetilde{\be}^U &=\bP_e^{-1}(\be^U - \bq_e)\\
\widetilde{\bh}^L &= \bP_h^{-1}(\bh^L - \bq_h), &\widetilde{\bh}^U &=\bP_h^{-1}(\bh^U - \bq_h) \label{eq:path-tilde-bounds}
\end{align}
Equations \eqref{eq:J-tilde}-\eqref{eq:path-tilde-bounds} constitute Problem $\widetilde{B}$.

\section{Necessary Conditions for the Scaled Problem}
Let $\blamt, \bmut$ and $\bnut$ be the adjoint, path and endpoint covectors respectively, associated with the necessary conditions for Problem $\widetilde{B}$. Then, it follows that the Hamiltonian, the Lagrangian of the Hamiltonian, and the Endpoint Lagrangian for Problem $\widetilde{B}$ are given by\cite{ross-book}:
\begin{itemize}
\item The Hamiltonian, $\widetilde{H}$:
\begin{align}\label{eq:H-tilde}
\widetilde{H}(\blamt, \bxt, \but, \widetilde{t}):= &\frac{p_t}{p_J} F(\bP_x \bxt(t) + \bq_x,\ \bP_u \but(t) + \bq_u,\ p_t \widetilde{t} +q_t) \nonumber \\
& \qquad + p_t \blamt^T\bP^{-1}_x \bff(\bP_x \bxt(t) + \bq_x,\ \bP_u \but(t) + \bq_u,\ p_t \widetilde{t} +q_t)
\end{align}
\item The Lagrangian of the Hamiltonian, $\overline{\widetilde{H}}$:
\begin{align}
\overline{\widetilde{H}}(\bmut, \blamt, \bxt, \but, \widetilde{t}) & :=  \widetilde{H}(\blamt, \bxt, \but, \widetilde{t}) \nonumber\\
 &\qquad + \bmut^T\bP_h^{-1}\big[\bh \big(\bP_x \bxt(t) + \bq_x,\ \bP_u \but(t) + \bq_u,\ p_t \widetilde{t} +q_t\big) - \bq_h \big] \label{eq:Bt-HL}
\end{align}

\item The Endpoint Lagrangian $\overline{\widetilde{E}}$:
\begin{align}
\overline{\widetilde{E}}(\bnut, \bxt_0, \bxt_f, \widetilde{t}_0, \widetilde{t}_f):=  \frac{1}{p_J} E(\bP_x \bxt_0 + \bq_x,\ \bP_x \bxt_f + \bq_x,\ p_t \widetilde{t}_0 +q_t,\ p_t \widetilde{t}_f + q_t) \nonumber\\
\qquad + \bnut^T\bP_e^{-1}\left[\be\big(\bP_x \bxt_0 + \bq_x,\ \bP_x \bxt_f + \bq_x,\ p_t \widetilde{t}_0 +q_t,\ p_t \widetilde{t}_f + q_t  \big) - \bq_e\right] \label{eq:Bt-Ebar}
\end{align}

\end{itemize}
The adjoint equation for Problem $\widetilde{B}$ is given by,
\begin{equation}
-\frac{d\blamt}{d\widetilde{t}} = \frac{\partial \overline{\widetilde{H}}}{\partial \bxt} \label{eq:Bt-adjoint}
\end{equation}
Evaluating the right-hand-side of \eqref{eq:Bt-adjoint} using \eqref{eq:Bt-HL} we get,
\begin{align}
-\frac{d\blamt}{d\widetilde{t}} = \frac{p_t}{p_J} \bP_x^T
\left[\frac{\partial F}{\partial\bx} \right] + p_t \left[\bP_x^{-1} \frac{\partial\bff}{\partial\bx}
\bP_x \right]^T \blamt
+ \left[\bP_h^{-1} \frac{\partial\bh}{\partial\bx}
\bP_x \right]^T \bmut \label{eq:Bt-ad-derive-1}
\end{align}

Following the same process for the stationarity condition associated with the Hamiltonian minimization condition for Problem $\widetilde{B}$, we get,
\begin{equation}
\frac{\partial \overline{\widetilde{H}}}{\partial\but} =  \frac{p_t}{p_J} \bP_u^T \left(\frac{\partial F}{\partial\bu}\right)
+ p_t \bP_u^T \left( \frac{\partial\bff}{\partial\bu} \right)^T \left(\bP_x^{-1}\right)^T \blamt
+ \bP_u^T \left( \frac{\partial\bh}{\partial\bu} \right)^T \left(\bP_h^{-1}\right)^T \bmut = \bzero
\end{equation}
Likewise the initial and final transversality conditions are given by,
\begin{subequations}
\begin{align}
-\blamt(\widetilde{t}_0) = \frac{\partial \overline{\widetilde{E}}}{\partial\bxt_0}= \frac{\bP_x^T}{p_J} \frac{\partial E}{\partial\bx_0} + \left[\bP_e^{-1}\frac{\partial\be}{\partial\bx_0}\bP_x\right]^T \bnut \label{eq:Bt-tvc-deriv-0}\\
\blamt(\widetilde{t}_f) = \frac{\partial \overline{\widetilde{E}}}{\partial\bxt_f}= \frac{\bP_x^T}{p_J} \frac{\partial E}{\partial\bx_f} + \left[\bP_e^{-1}\frac{\partial\be}{\partial\bx_f}\bP_x\right]^T \bnut \label{eq:Bt-tvc-deriv-f}
\end{align}
\end{subequations}
Finally, the complementarity conditions are given by $\bmut\dagger\bht$ and $\bnut\dagger\bet$.

\section{The Balancing Equations}

\emph{\textbf{Proposition A:}}
Let, $\bxt^*(\cdot), \but^*(\cdot), \widetilde{t}^*_0$ and $\widetilde{t}^*_f$ be an extremal solution to the scaled problem $\widetilde{B}$.  Let, $\blamt^*(\cdot), \bmut^*(\cdot)$ and $\bnut^*$ be a multiplier triple associated with this extremal solution.  Let $\bP_\lambda$, $\bP_\mu$ and $\bP_\nu$ be invertible matrices defined according to:
\begin{subequations}\label{eq:result-main-def}
    \begin{align}
    \bP_\lambda &:= p_J \left[\bP_x^{-1}\right]^T \\[1em]
    \bP_\mu &:= \frac{p_J}{p_t} \left[\bP_h^{-1}\right]^T \\[0.5em]
    \bP_\nu &:= p_J \left[\bP_e^{-1}\right]^T
    \end{align}
\end{subequations}
Then an extremal solution to the unscaled Problem $B$ exists and is given by:
\begin{subequations}\label{eq:result-main}
\begin{align}
\bx^*(\cdot) := & \bP_x\ \bxt^*(\cdot) + \bq_x \label{eq:result-primal-x}\\
\bu^*(\cdot) := & \bP_u\ \but^*(\cdot) + \bq_u \\
t := & p_t\ \widetilde{t} + q_t \label{eq:result-primal-t}\\
\blam^*(\cdot) := & \bP_\lambda\ \tlam^*(\cdot) \label{eq:result-dual-lam}\\
\bmu^*(\cdot) := & \bP_\mu\ \tmu^*(\cdot)\\
\bnu^* := & \bP_\nu\ \tnu^* \label{eq:result-dual-nu}
\end{align}
\end{subequations}
%
\emph{\textbf{Proof:}}
The proof of \eqref{eq:result-primal-x}--\eqref{eq:result-primal-t} follows quite simply by construction.  The proof \eqref{eq:result-dual-lam}--\eqref{eq:result-dual-nu} follows by substituting \eqref{eq:result-main} in Problem $B^\lambda$ and using \eqref{eq:result-main-def} to show that the resulting equations are the same as the necessary conditions for the scaled problem derived in the previous section.

\subsection*{Remark 1}
Once $\bP_x, \bP_h $ and $\bP_e$ are chosen to scale the primal variables and constraints, then, according to Proposition~$A$,  there exists dual variables $\blam^*(\cdot), \bmu^*(\cdot)$, and $\bnu^*$, that get scaled automatically in compliance with \eqref{eq:result-main}.  Consequently balancing can now be defined more precisely as choosing $\bP_x, \bP_h $ and $\bP_e$ along with $p_J$ and $p_t$ such that the values of the covectors $\blam^*(\cdot), \bmu^*(\cdot)$, and $\bnu^*$ are of similar orders of magnitude as their corresponding vectors.  In contrast, scaling is choosing $\bP_x, \bP_h $ and $\bP_e$ such that the values of the corresponding vectors (as well as their components) are of similar magnitude relative to each other.  When both requirements are met, the problem is said to be scaled and balanced.

\subsection*{Remark 2}
Substituting \eqref{eq:result-main-def} in \eqref{eq:H-tilde}, it is clear that the value of the Hamiltonian transforms according to
\begin{equation}\label{eq:HVT}
H (\blam, \bx, \bu, t) = \left(\frac{p_J}{p_t}\right)\widetilde{H} (\blamt, \bxt, \but, \widetilde{t})
\end{equation}

\subsection*{Remark 3}
A natural choice for  $\bP_x, \bP_h $ and $\bP_e$ are diagonal matrices. For these choices, $\bP_\lambda, \bP_\mu $ and $\bP_\nu$ are also diagonal matrices. In this situation, each component of a vector is independently related to each component of its corresponding covector. This fact can be utilized in a numerical setting for a simple algorithmic technique for scaling and balancing.

\subsection*{Remark 4}
Let the units of $\bx, \bh$ and $\be$ be given according to,
\begin{subequations}\label{eq:vec-def-units}
\begin{align}
\renewcommand\arraystretch{1.0} \bx &:= \left[
           \begin{array}{c}
             x_1 \\
             x_2 \\
             \vdots \\
             x_{N_x}
           \end{array}
         \right]   \ \begin{array}{c}
                   x_1 \text{-units} \\
                   x_2 \text{-units} \\
                   \vdots \\
                   x_{N_x} \text{-units}
                 \end{array}\label{eq:state-def-units}
                 \\[1em]
\bh &:= \left[
           \begin{array}{c}
             h_1(\bx, \bu, t) \\[0.7em]
             h_2(\bx, \bu, t) \\ 
             \vdots \\
             h_{N_h}(\bx, \bu, t)
           \end{array}
         \right]   \ \begin{array}{c}
                   h_1 \text{-units} \\[1em]
                   h_2 \text{-units} \\ 
                   \vdots \\
                   h_{N_h} \text{-units}
                 \end{array}\label{eq:path-def-units}
\\[1em]
\be(\bx_0, \bx_f, t_0, t_f) &:= \left[
           \begin{array}{c}
             e_1(\bx_0, \bx_f, t_0, t_f) \\
             e_2(\bx_0, \bx_f, t_0, t_f) \\
             \vdots \\
             e_{N_e}(\bx_0, \bx_f, t_0, t_f)
           \end{array}
         \right]   \ \begin{array}{c}
                   e_1 \text{-units} \\
                   e_2 \text{-units} \\
                   \vdots \\
                   e_{N_e} \text{-units}
                 \end{array}
\end{align}
\end{subequations}
Note that no assumption is made on any of the units in \eqref{eq:vec-def-units}.  Thus, for example, $x_1$-units may be meters and $x_2$-units may be feet even if the pair $(x_1, x_2)$ is the position coordinate of the same point mass.  Obviously, the Euclidian norm of the numbers given by $x_1$ and $x_2$ has no physical meaning.  In the same spirit, $x_3$-units may be yards per day even if the variable $x_3$ is the time rate change of $x_1$.  Despite such arbitrary choices, it is still possible to measure a ``length'' of a vector through the use of a covector.  To this end, we let $\bP_x, \bP_h $ and $\bP_e$ be diagonal matrices.  Then, it follows from \eqref{eq:result-main-def} and \eqref{eq:result-main} that the covectors have units given by,
\begin{subequations}\label{eq:covec-def-units}
\begin{align}\renewcommand\arraystretch{1.0} \blam &:= \left[
           \begin{array}{c}
             \lambda_1 \\
             \lambda_2 \\
             \vdots \\
             \lambda_{N_x}
           \end{array}
         \right]   \ \begin{array}{c}
                   CU/x_1 \text{-units} \\
                   CU/x_2 \text{-units} \\
                   \vdots \\
                   CU/x_{N_x} \text{-units}
                 \end{array}\label{eq:costate-def-units}
\\
\bmu &:= \left[
           \begin{array}{c}
             \mu_1 \\[0.7em]
             \mu_2 \\ 
             \vdots \\
             \mu_{N_h}
           \end{array}
         \right]   \ \begin{array}{c}
                   \frac{CU/TU}{h_1 \text{-units}} \\[1em]
                   \frac{CU/TU}{h_2 \text{-units}} \\ 
                   \vdots \\
                   \frac{CU/TU}{h_{N_h} \text{-units}}
                 \end{array}\label{eq:mu-covector-def-units}
\\
\bnu &:= \left[
           \begin{array}{c}
             \nu_1 \\
             \nu_2 \\
             \vdots \\
             \nu_{N_e}
           \end{array}
         \right]   \ \begin{array}{c}
                   CU/e_1 \text{-units} \\
                   CU/e_2 \text{-units} \\
                   \vdots \\
                   CU/e_{N_e} \text{-units}
                 \end{array}\label{eq:nu-covector-def-units}
\end{align}
\end{subequations}
where, $CU$ is the cost unit, and $TU$ is the time unit. Equation \eqref{eq:covec-def-units} was first introduced in \cite{ross-book} as part of the definition of a covector associated with the relevant vector.  Note also from \eqref{eq:covec-def-units} that a covector always has some unit of frequency; i.e., it is always given in terms of some common unit per some unit.  The common unit in \eqref{eq:costate-def-units} and \eqref{eq:nu-covector-def-units} is the cost unit, but the common unit in \eqref{eq:mu-covector-def-units} is the cost unit per time unit.  Consequently, we can now measure the ``length'' of a vector by an appropriate covector.  For instance, $\bx^T\blam$ generates a scalar in terms of $CU$'s despite that its Euclidean norm might not be computable due to the disparity in units of the constituents of $\bx$. See \cite{ross-book}, Sec.~2.2 for further details.

\subsection*{Remark 5}

As a consequence of \eqref{eq:vec-def-units} and \eqref{eq:covec-def-units}, scaling may be conceived as simply changing units. Consequently, ``nondimensionalization'' is also changing units, and should not be construed as eliminating units.

\subsection*{Remark 6}

\emph{\textbf{The Hamiltonian is not dimensionless. The unit of the Hamiltonian is cost unit per time unit}}.  This follows from \eqref{eq:covec-def-units}.  When the cost unit is the same as the time unit (e.g. time optimality), then, and only then, is the Hamiltonian truly dimensionless.

\subsection*{Remark 7}
Using \eqref{eq:dynamics-tilde} we can write,
\begin{equation}\label{eq:dyn-jac-X}
\frac{\partial\widetilde{\bff}}{\partial \bxt} = p_t\bP_x^{-1} \left(\frac{\partial\bff}{\partial \bx}\right)\bP_x
\end{equation}
Because of the similarity transformation on the right-hand side of \eqref{eq:dyn-jac-X}, the spectral radii of the Jacobians are related by,
\begin{equation}\label{eq:spectral-radii-X}
\rho\left(\partial_{\widetilde{x}}\widetilde{\bff}\right) = p_t \rho\big( \partial_x\bff\big)
\end{equation}
where, $\rho(\cdot)$ is the spectral radius of $(\cdot)$.
Rewriting $p_t$ as $(t_f-t_0)/(\widetilde{t}_f-\widetilde{t}_0)$ and substituting in \eqref{eq:spectral-radii-X}, we get an invariance equation,
\begin{equation}\label{eq:sensitivity-invar}
\left(\widetilde{t}_f-\widetilde{t}_0\right)\rho\left(\partial_{\widetilde{x}}\widetilde{\bff}\right) = \big(t_f-t_0\big) \rho\big(\partial_x\bff\big)
\end{equation}
Because the product of the spectral radius and the time horizon is a key sensitivity factor (see Ref.~\cite{ross-book}, Sec.~2.9), it follows from \eqref{eq:sensitivity-invar} that \emph{\textbf{the curse of sensitivity cannot be mitigated by scaling}}.

\section{Example  Illustrating Designer Units, Scaling and Balancing}
\begin{figure}[h!]
      \centering
      {\includegraphics[clip,width = 0.55\textwidth]{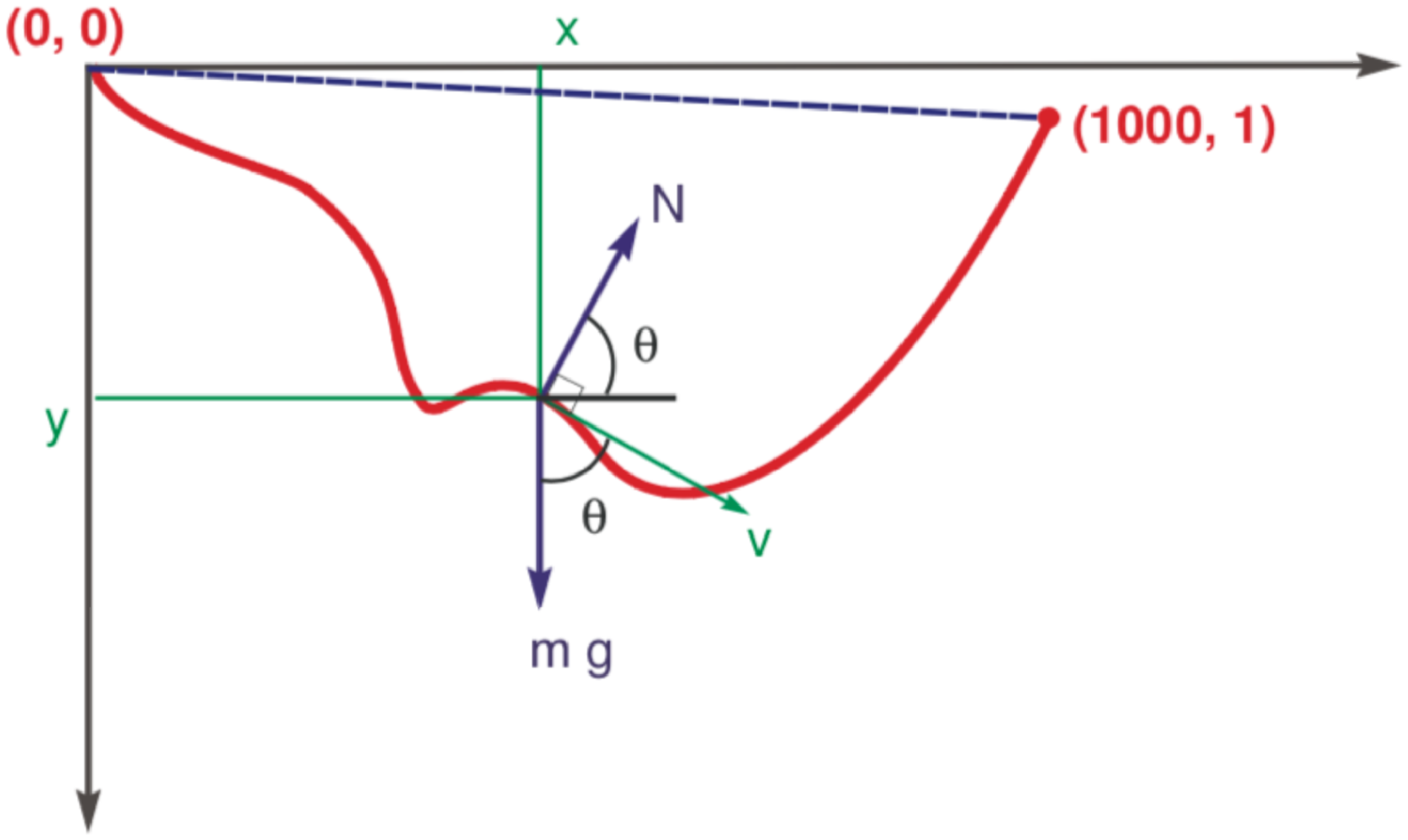}}
      \caption{\textsf{Schematic for the ``bad'' Brachistochrone problem.  Figure is not to scale.}}
    \label{fig:BB-sketch}
   \end{figure}

Although the concepts of designer units, scaling and balancing  have been used for more than a decade to generate successful flight and field operations\cite{bhatt:opm,SIAMnews,zpm:NASA-report,zpm:IEEE, TRACE-IEEE-Spectrum,CDC-Workshop,PSReview-ARC-2012}, the process has largely been ad hoc until recently. Because the flight problems are too involved for illustrating the procedures, we design a far simpler problem that typifies the systematic process of scaling and balancing. To this end, consider the well-known Brachistochrone problem. Instead of using well-scaled numbers as is customarily discussed in many textbooks, we purposefully choose the final-time condition on the position coordinates $(x, y)$ to be widely disparate and given by $(1000, 1)$ meters; see Fig.~\ref{fig:BB-sketch}.
This results in a badly-scaled or ``bad'' Brachistochrone problem that can be formulated as\cite{ross-book},
\[
\left.\begin{aligned}
\hspace{0.8cm}\phantom{\quad \textsf{(preamble)}}
\X &=\real{3}   &\U &= \Real& \\
\bx &= (x, y, v)  \quad &\bu &= \theta & \\
\end{aligned}\hspace{1.2cm} \right\} \hspace{0.3cm} \text{(preamble)}
\]
\[
\begin{aligned}
\renewcommand*{\leftterm}{(x_f - x^f, y_f -y^f)}
\renewcommand*{\rightterm}{\phantom{(0, 0, 0, 0, 0)}}
\settowidth{\leftside}{\term{\leftterm}}
\settowidth{\rightside}{\term{\rightterm}}
(B_R) \left\{
\begin{array}{ll}
\text{Minimize }&
 \left.\begin{aligned}
 \makebox[\leftside][r]{\term{J[\bx(\cdot), \bu(\cdot), t_f]}} & = \makebox[\rightside][l]{\term{t_f}}\\
 \end{aligned}\right\} \quad\, \text{(cost)}\\[0.7em]
\text{Subject to} &  
\\[-1.5em]
&  \left.\begin{aligned}
    \makebox[\leftside][r]{\term{\dot x}} &= \makebox[\rightside][l]{\term{v\sin\theta}}\\
    \dot y &= v\cos\theta \\
    \dot v &= g\cos\theta
  \end{aligned}\right\} \quad \text{(dynamics)} \\
&  \left.\begin{aligned}
    \makebox[\leftside][r]{\term{(t_0, x_0, y_0, v_0)}} &= \makebox[\rightside][l]{\term{(0, 0, 0, 0)}} \\
    (x_f, y_f) &=  (1000, 1)
  \end{aligned}\right\} \quad \text{(endpoints)}
  %
\end{array}
\right.
\end{aligned}
\]
where, $g = 9.8\ \textit{meters}/\textit{sec}^2$.  See Fig.~\ref{fig:BB-sketch} for a physical definition of the variables.

%
%

\subsection{Illustrating the Process for Choosing Designer Units}
As a consequence of \eqref{eq:state-def-units}, we can write,
\begin{equation}\renewcommand\arraystretch{1.0} \bx := \left[
           \begin{array}{c}
             x \\
             y \\
             v
           \end{array}
         \right]   \ \begin{array}{c}
                   x\text{-units} \\
                   y\text{-units} \\
                   v\text{-units}
                 \end{array}\label{eq:state-def-units-Brac1}
\end{equation}
Because the cost unit is the same as the time unit, the adjoint covector is defined by
\begin{equation}\renewcommand\arraystretch{1.0} \blam := \left[
           \begin{array}{c}
             \lambda_x \\
             \lambda_y \\
             \lambda_{v}
           \end{array}
         \right]   \ \begin{array}{c}
                   t\text{-units}/x\text{-units} \\
                   t\text{-units}/y\text{-units} \\
                   t\text{-units}/v \text{-units}
                 \end{array}\label{eq:costate-def-brac:1}
\end{equation}
As a result, the Hamiltonian
\begin{equation}\label{eq:Bb-Hamil}
H(\blam, \bx, \bu) = \lambda_x v \sin\theta + \lambda_y v \cos\theta + \lambda_v g \cos\theta
\end{equation}
is dimensionless.

Applying the Hamiltonian minimization condition, we get
\begin{equation}\label{eq:Bb-HMC}
\lambda_x(t) v(t) \cos\theta(t) - \lambda_y(t) v(t) \sin\theta(t) - \lambda_v(t) g \sin\theta(t) = 0 \quad( \forall\ t \in [t_0, t_f] )
\end{equation}
The adjoint equations,
\begin{equation}\label{eq:Bb-adjoint}
\begin{aligned}
\dot\lambda_x = & 0 \\
\dot\lambda_y = & 0 \\
\dot\lambda_v = & -\lambda_x \sin\theta - \lambda_y\cos\theta
\end{aligned}
\end{equation}
indicate that $\lambda_x$ and $\lambda_y$ are constants. In addition, the transversality condition,
\begin{equation}
\lambda_v(t_f) = 0
\end{equation}
and the Hamiltonian value condition,
\begin{equation}\label{eq:Bb-HVC-f}
H[@t_f] := \lambda_x(t_f) v(t_f) \sin\theta(t_f) + \lambda_y(t_f) v(t_f) \cos\theta(t_f) + \lambda_v(t_f) g \cos\theta(t_f) = -1
\end{equation}
complete the computational set of conditions that define the boundary value problem.
These equations can also be used as part of the totality of a verification and validation  of a candidate optimal solution obtained by any computational method.

In order to perform initial scaling, we simply take the given numerical data as a starting point.   For the numerics given in Problem $B_R$,  we expect $x$ to satisfy,
\begin{equation}\label{eq:search-x}
 0 \le x \le 1000
\end{equation}
In the absence of further analysis, we can assume $y$ to take on similar range of values. Furthermore, because the $x$-distance is relatively large, we expect the time of travel to be relatively large (in terms of seconds).  Based on these heuristics, we choose an initial set of scaling factors according to:
\begin{equation}\label{eq:BB-designer-units-1}
\begin{aligned}
 x &=  P_x\ \widetilde{x} &&= 100\ \widetilde{x}  \\
 y & = P_y\ \widetilde{y} &&= 20\ \widetilde{y}  \\
 v & = P_v\ \widetilde{v} &&= 10\ \widetilde{v}  \\
 \theta & = P_\theta\ \widetilde{\theta} && = \widetilde{\theta}  \\
 t & = p_t\ \widetilde{t} && = 10\ \widetilde{t}
\end{aligned}
\end{equation}
As a result of \eqref{eq:BB-designer-units-1}, the large variation in $x$ indicated by \eqref{eq:search-x} is tempered by $\widetilde{x}$ according to
$$0 \le \widetilde{x} \le 10$$
Note that the numbers given in \eqref{eq:BB-designer-units-1} imply a unit of distance along the $x$-axis that is completely different from the unit of distance along the $y$ axis! In fact, these numbers constitute a specific system of units that do not conform with the metric system or any other set of the standard units; hence, these are designer units. The conversion between the designer units of \eqref{eq:BB-designer-units-1} and the metric units is given by,
\begin{equation}
\begin{aligned}
 1\ \widetilde{x}\text{-unit} &= \text{1 unit of distance along $x$-axis}  &&= 100 \text{ meters} \\
 1\ \widetilde{y}\text{-unit} &= \text{1 unit of distance along $y$-axis}  &&= 20 \text{ meters} \\
 1\ \widetilde{v}\text{-unit} &= \text{1 unit of speed} && = 10 \text{ meters/second} \\
 1\ \widetilde{t}\text{-unit} & =\text{1 unit of time } && = 10  \text{ seconds}
\end{aligned}
\end{equation}
It is important to note that the velocity unit is completely independent of any of the $x$, $y$ or $t$ units.  Consequently, these designer units are not consistent in the sense that,
$$\frac{d\widetilde{y}}{d\widetilde{t}} \ne \widetilde{v} \cos\widetilde{\theta} $$
To drive home this point, we note that we can no longer express $g$ in terms of distance units per the square of time units.  For instance, in the metric system, the unit of $g$ is given by $\textit{meters}/\textit{sec}^2$.  Because we chose distance units along the $x$- and $y$-directions to be independent of each other, it is clear that that we cannot regard $g$ in terms of ``distance units per the square of time units.''  The proper unit for $g$ is obtained by considering $d\widetilde{v}/d\widetilde{t}$. This implies that we may regard the gravitational acceleration as being transformed according to
\begin{equation}\label{eq:g-units}
\widetilde{g} = \left(\frac{p_t}{P_v}\right) g  = 9.8\ \widetilde{v}\text{-unit}/\widetilde{t}\text{-unit}
\end{equation}
That is, $\widetilde{g} \ne 1$ numerically, which is a typical number for canonical units\cite{longuski}.  That $\widetilde{g} = g $ numerically in \eqref{eq:g-units} is simply coincidental and as a result of choosing $p_t = P_v$ in \eqref{eq:BB-designer-units-1}.

Using the scaling units of \eqref{eq:BB-designer-units-1} the endpoint conditions can be written as,
\begin{equation}\label{eq:Bb-xy-final-t}
\begin{aligned}
\widetilde{t}_0 &= 0    && &&(\text{in }\widetilde{t}\text{-units})\\
\widetilde{x}(\widetilde{t}_0) &= 0 &\widetilde{x}(\widetilde{t}_f) &= 10   &&(\text{in }\widetilde{x}\text{-units})  \\
\widetilde{y}(\widetilde{t}_0) &= 0 &\widetilde{y}(\widetilde{t}_f) &= 0.05  \qquad &&(\text{in }\widetilde{y}\text{-units}) \\
\widetilde{v}_0 &= 0    && &&(\text{in }\widetilde{v}\text{-units})
\end{aligned}
\end{equation}
Imposing the endpoint conditions according to \eqref{eq:Bb-xy-final-t} is tantamount to choosing $\bP_e$ according to,
\begin{equation}\label{eq:Bb-Pe}
\bP_e = \textrm{diag}(10, 100, 100, 20, 20, 10)
\end{equation}
This follows as a direct consequence of \eqref{eq:BB-designer-units-1}.
Finally, we scale the cost functional using $p_J = 10$ so that we can write,
$$\widetilde{J}[\bxt(\cdot), \but(\cdot), \widetilde{t}_0, \widetilde{t}_f ] =  \widetilde{t}_f$$

At this juncture, we wish to emphasize that the purpose of initial scaling is not necessarily to generate the best set of designer units; rather, it is largely directed at producing a work-flow for balancing.  As will be apparent shortly, once an initial numerical result is obtained, balancing can usually be performed in just about two iterations. One situation when initial scaling becomes critical to the work-flow is when no numerical result is achieved simply because of poor scaling.  The detection and mitigation of this problem are open areas of research.

\subsection{Illustrating the Process of Descaling the Covectors }
For the purposes of clarity of the discussions to follow, we use the following terminology:
\begin{enumerate}
\item Unscaled: This refers to all numbers and variables associated with the original (or unscaled) problem.
\item Scaled: This refers to all numbers and variables associated with the affinely transformed (or scaled) problem.
\item Descaled: This refers to all numbers and variables that are purported solutions to the unscaled problem obtained via \eqref{eq:result-main} and a solution to the scaled problem.
\end{enumerate}
Applying \eqref{eq:result-main} to descale the adjoint covector, we get,
\begin{equation}\label{eq:Bb-lam-descale}
\begin{aligned}
\lambda_x &= \widetilde{\lambda}_x \left(\frac{p_J}{P_x}\right) &=&& \frac{\widetilde{\lambda}_x}{10} &&& \left(\frac{\text{seconds}}{\text{meters}} \right)\\
\lambda_y &= \widetilde{\lambda}_y \left(\frac{p_J}{P_y}\right)  &=&& \frac{\widetilde{\lambda}_y}{2} &&& \left(\frac{\text{seconds}}{\text{meters}} \right)\\
\lambda_v &= \widetilde{\lambda}_v \left(\frac{p_J}{P_v}\right)  &=&& \widetilde{\lambda}_v &&& \left(\frac{\text{seconds}^2}{\text{meters}} \right)
\end{aligned}
\end{equation}
Similarly, from \eqref{eq:Bb-Pe} and \eqref{eq:result-main-def} we have
\begin{equation}\label{eq:BB-nu-descale}
\bP_\nu = 10 \left[\bP_e^{-1}\right]^T  \Rightarrow (\nu_{t_0}, \nu_{x_0}, \nu_{x_f}, \nu_{y_0}, \nu_{y_f}, \nu_{v_0}) = \left(\widetilde{\nu}_{t_0}, \frac{\widetilde{\nu}_{x_0}}{10}, \frac{\widetilde{\nu}_{x_f}}{10}, \frac{\widetilde{\nu}_{y_0}}{2}, \frac{\widetilde{\nu}_{y_f}}{2}, \widetilde{\nu}_{v_0}\right)
\end{equation}
where, $\nu$ and $\widetilde{\nu}$ with the appropriate subscripts are the endpoint multipliers associated with the initial and final-time conditions.  Note that these endpoint multipliers also have units similar to those identified in \eqref{eq:Bb-lam-descale}.

\subsection{Illustrating the Numerical Process of Scaling and Balancing}
All of the analysis so far has been agnostic to the specific choice of a numerical method or software. To demonstrate the numerical process, any appropriate mathematical software may be used. We begin by choosing DIDO${}^\copyright$, a state-of-the-art MATLAB${}^\circledR$ toolbox for solving optimal control problems\cite{ross-book}. DIDO is the same tool that was used in all of the flight applications noted earlier. It is based on the the spectral algorithm\cite{knots,autoknots,kang-ijrc,spec-alg, kang-rate, arb-grid} for pseudospectral optimal control and does not require any guess of the solution to solve the problem. Furthermore, DIDO automatically generates all of the covectors associated with a generic optimal control problem (see Problem $B^\lambda$ presented in Section II of this paper) through an implementation of the Covector Mapping Principle\cite{ross-book,PSReview-ARC-2012,arb-grid}. That is, DIDO generates a guess-free candidate solution to the BVP while only requiring the data functions for Problem $B$.  This is why the spectral algorithm and its implementation in DIDO do not belong to the class of ``direct'' or ``indirect'' methods.  In fact, these ideas effectively obviate the need for such a classification; see Sec.~2.9.2 of Ref.~\cite{ross-book}.

\subsubsection{Initial Scaling and Descaling}
\begin{figure}[h!]
      \centering
      {\includegraphics[angle=0, width = 0.49\textwidth]{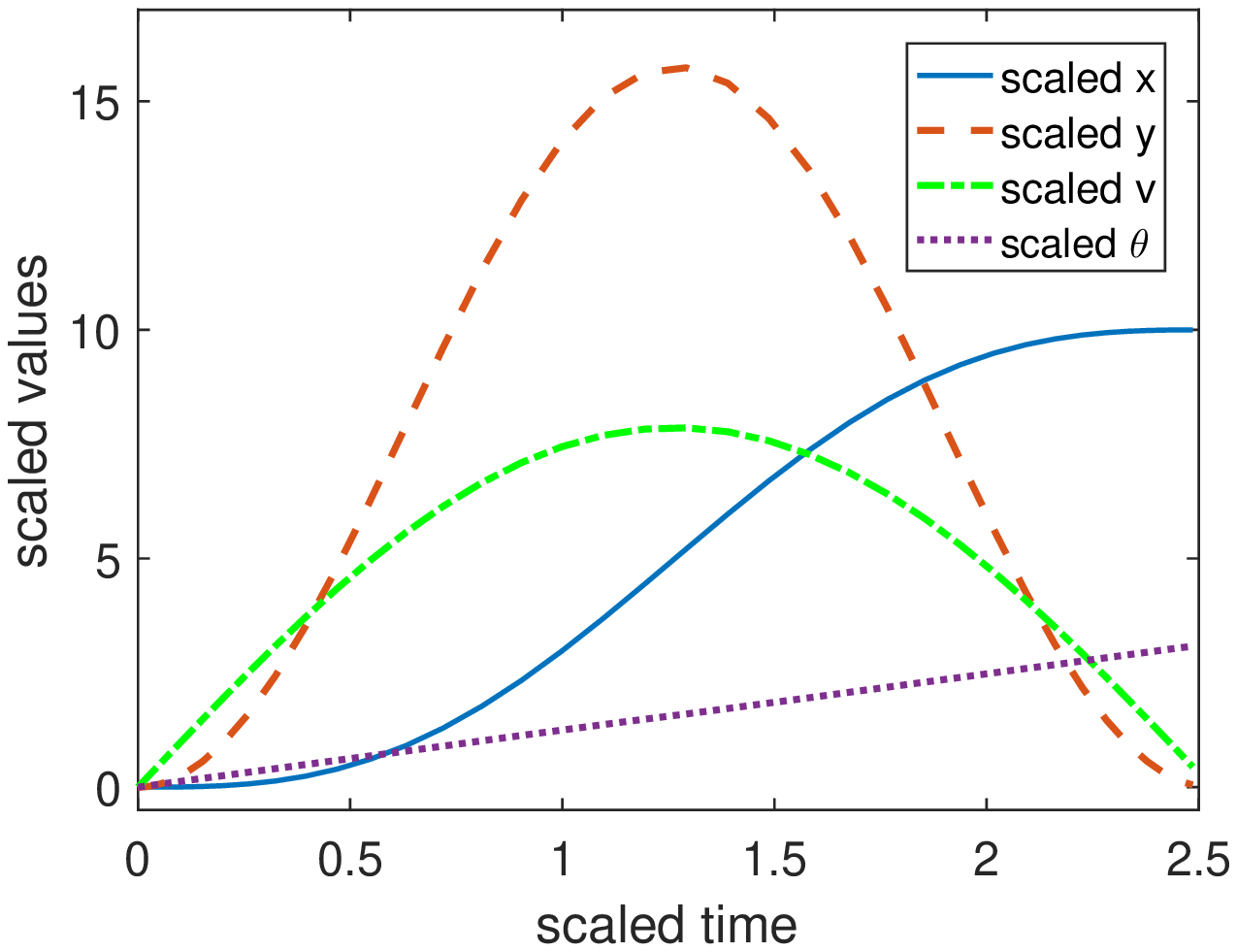}}
            {\includegraphics[angle=0, width = 0.49\textwidth]{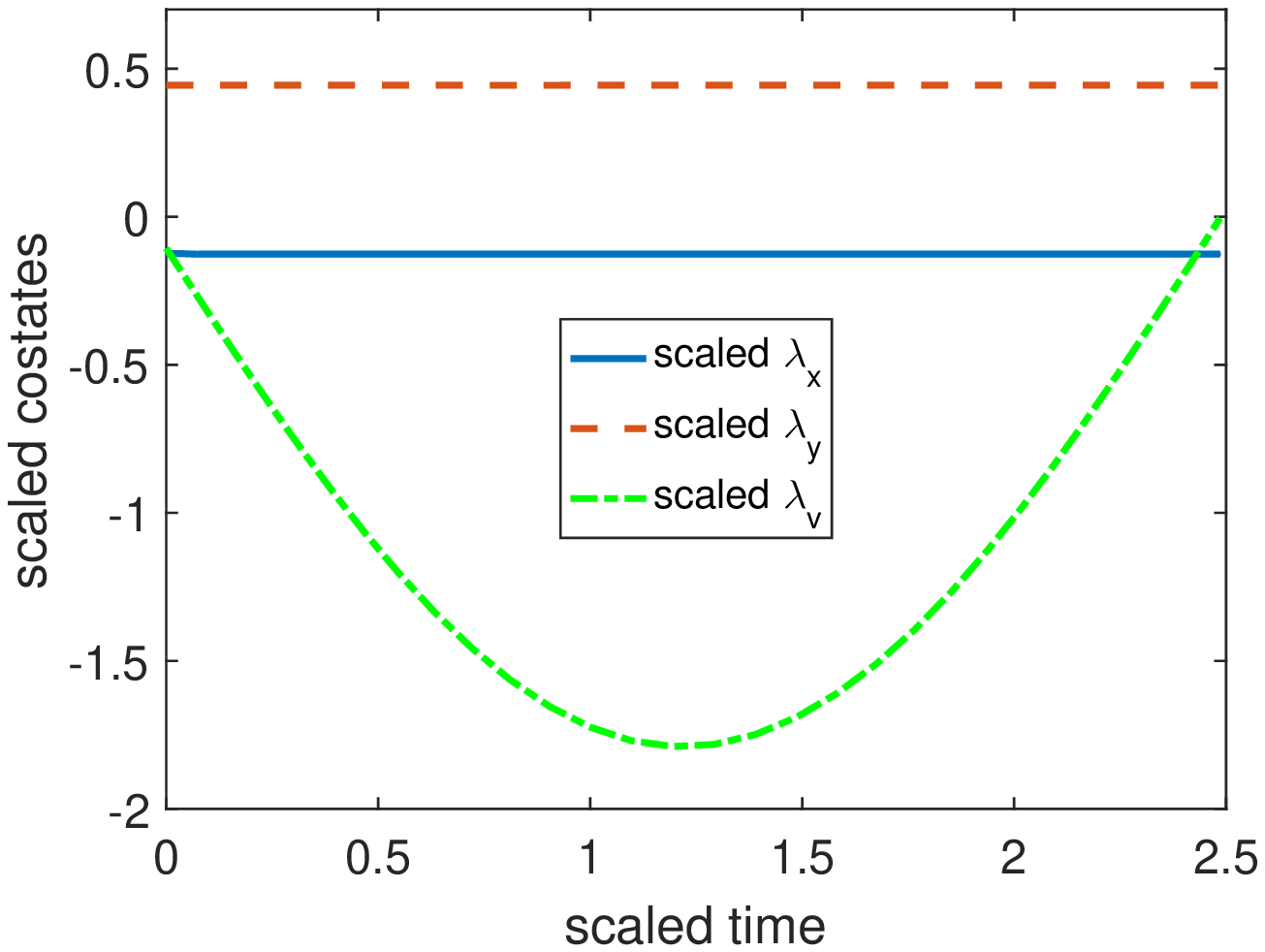}}
      \caption{\textsf{Guess-free primal (left) and dual (right) solutions to the badly-scaled Brachistochrone problem.}}
    \label{fig:BB-scaledPD}
   \end{figure}

\begin{figure}[h!]
      \centering
      {\includegraphics[angle=0, width = 0.47\textwidth]{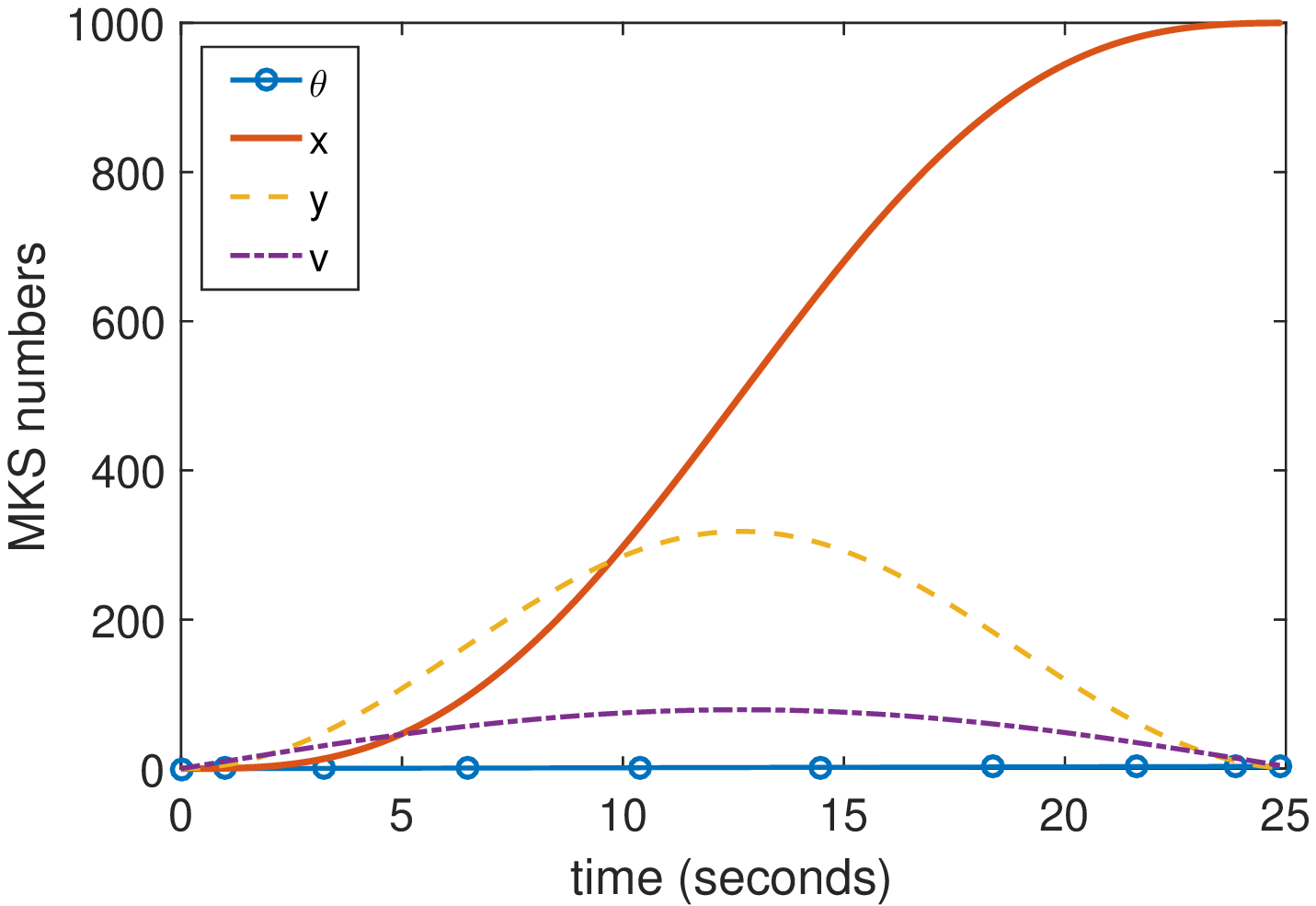}}
            {\includegraphics[angle=0, width = 0.5\textwidth]{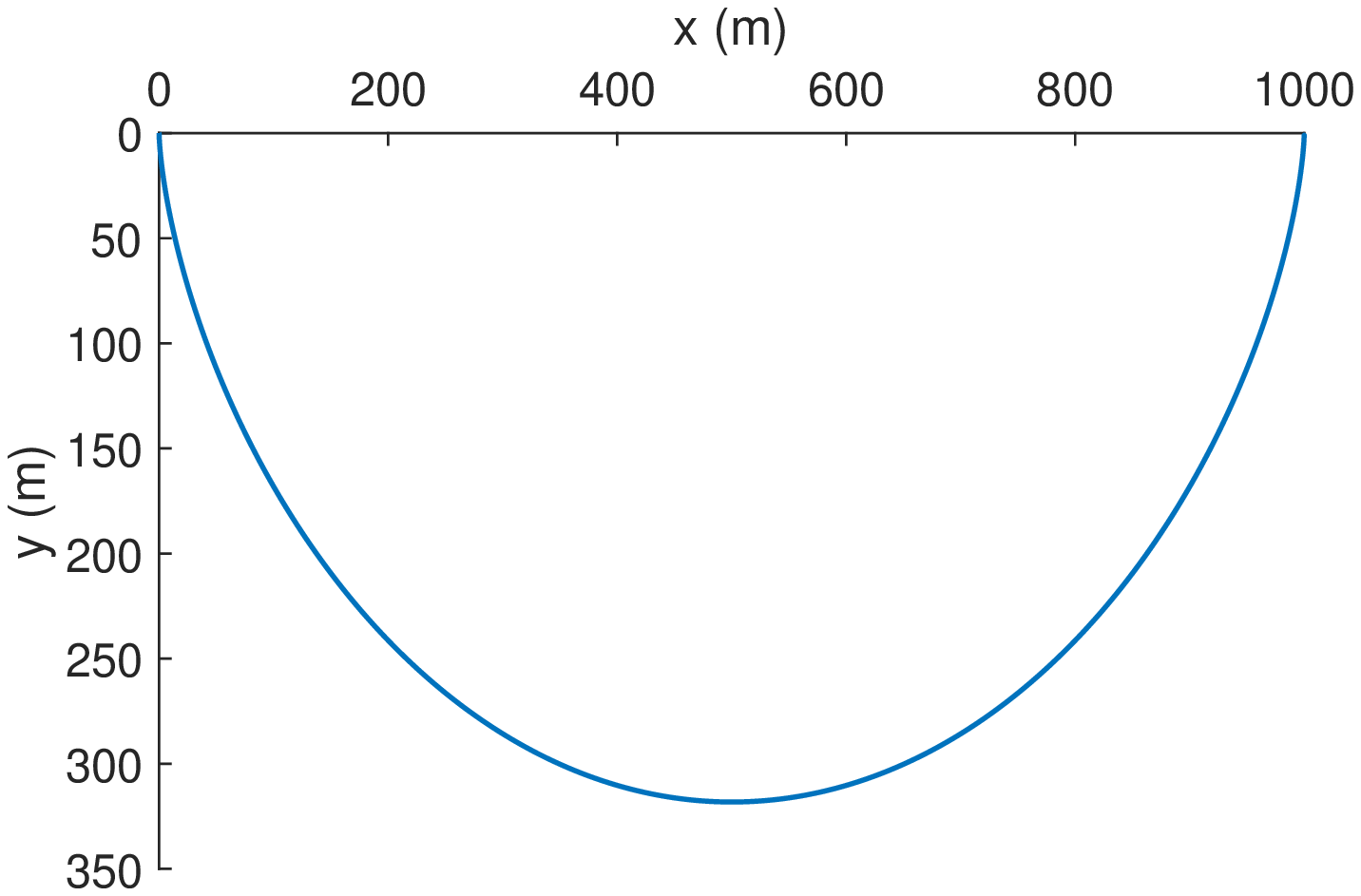}}
      \caption{\textsf{Primal feasible solution to the badly-scaled Brachistochrone problem indicating variations in three-orders of magnitude.}}
    \label{fig:BB-primalProp}
   \end{figure}

A candidate primal-dual solution (generated by DIDO) is shown in Fig.~\ref{fig:BB-scaledPD}.
%
%
By definition, these are simply candidate solutions to the scaled problem. They have not yet been validated.
Because it is frequently more meaningful to validate results in physical units, we first descale the primal variables using \eqref{eq:BB-designer-units-1}.  The descaled candidate control trajectory is then used to propagate the initial conditions,
$$ x(0) = 0, \quad y(0) = 0, \quad v(0) = 0  $$
using linear interpolation for the controls and \textsf{ode45} in MATLAB.  The propagated state trajectory is shown in Fig.~\ref{fig:BB-primalProp}.
%
%
The propagated values of $x(t_f)$ and $y(t_f)$ satisfy the final-time conditions to the following precision,
$$\mid x(t_f) - 1000\mid = 3.6\times 10^{-3} \ m , \quad \mid y(t_f)-1 \mid = 1.9 \times 10^{-4} \ m  $$
Thus, the descaled solution is verifiably feasible. In flight applications, such an independent verification of feasibility is critical to a successful pre-flight checkout\cite{zpm:IEEE,TRACE-IEEE-Spectrum,PSReview-ARC-2012}.

\begin{figure}[h!]
      \centering
      {\includegraphics[angle=0, width = 0.49\textwidth]{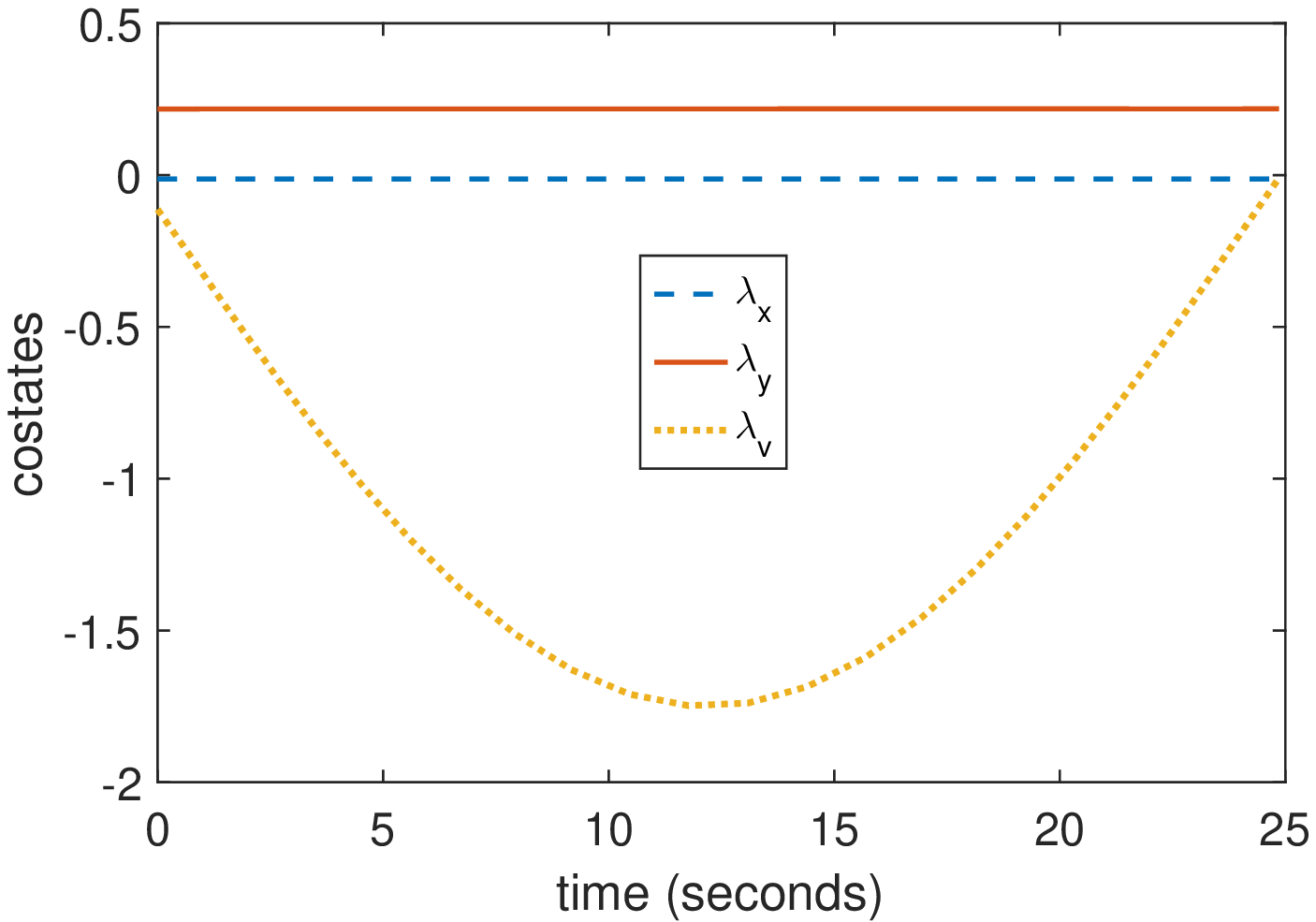}}
            {\includegraphics[angle=0, width = 0.49\textwidth]{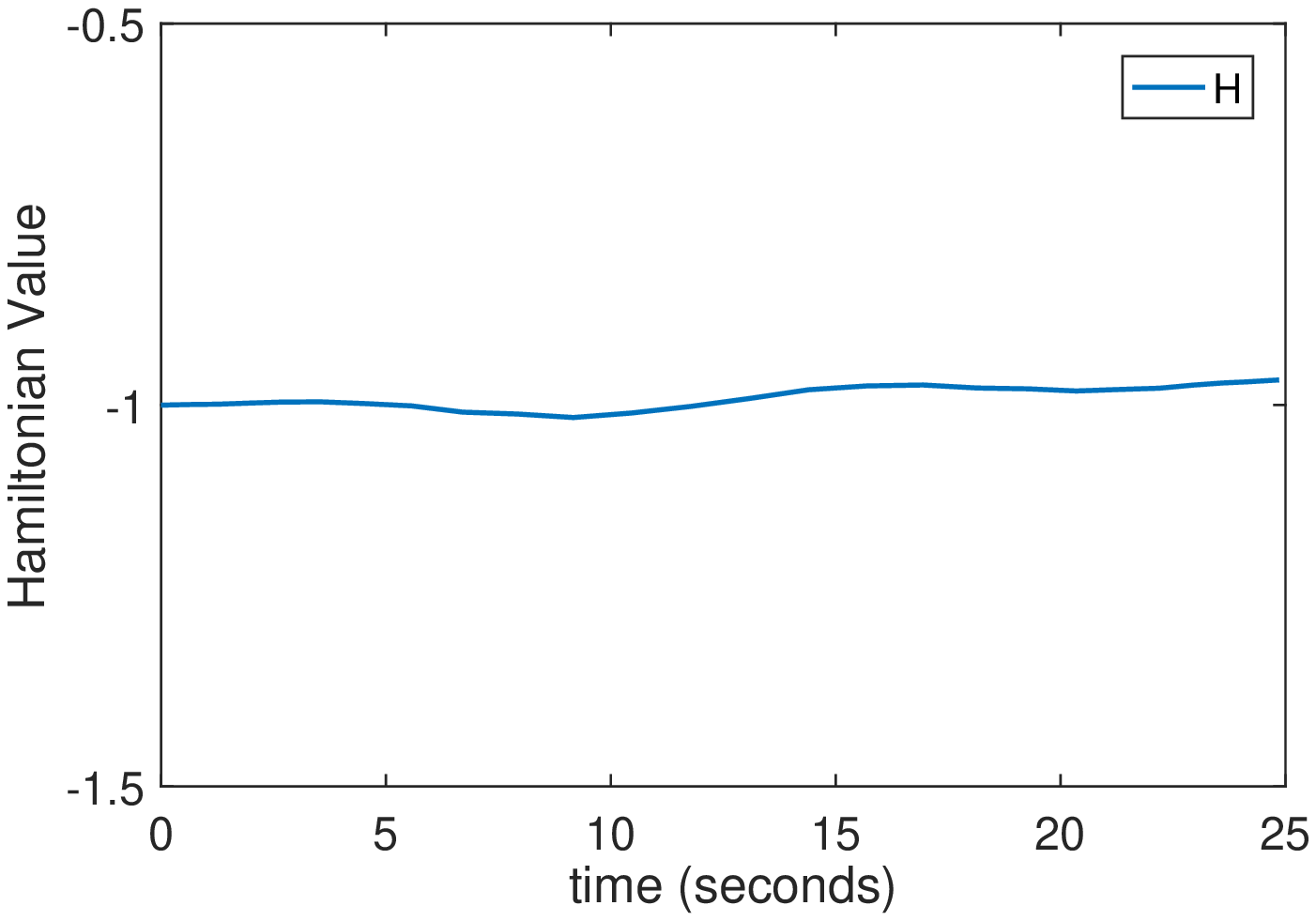}}
      \caption{\textsf{Descaled costates and the evolution of the Hamiltonian for the badly-scaled Brachistochrone problem.}}
    \label{fig:BB-costates}
   \end{figure}

To validate the extremality of the feasible solution, we use \eqref{eq:Bb-lam-descale} and \eqref{eq:HVT} to descale the adjoint covectors and the evolution of the Hamiltonian respectively.  The results are shown in Fig.~\ref{fig:BB-costates}.
%
%
It is apparent that $\lambda_x(t)$ and $\lambda_y(t)$ are constants as required by \eqref{eq:Bb-adjoint}.  It is also apparent that the Hamiltonian is nearly a constant and equal to $-1$ as required by \eqref{eq:Bb-HVC-f} and the first integral.
Thus, the theoretical necessary conditions are satisfied up to the indicated approximations.\footnote{Although the spectral algorithm can theoretically generate very accurate solutions\cite{ross-book,PSReview-ARC-2012,spec-alg,arb-grid}, the Hamiltonian evolution equation is satisfied only weakly\cite{lncis}; hence, the Hamiltonian is not expected to be equal to $-1$ in the strong $L_\infty$-norm in Fig.~\ref{fig:BB-costates}.   In addition, because no Jacobian information was provided in the generation of Fig.~\ref{fig:BB-costates}, the accuracy in the computation of the dual variables is expected to be lower than that of the primal solution. Note also that dual (or primal) tolerances cannot be set to arbitrarily small numbers (e.g., $10^{-6}$ or $10^{-8}$) to attain ``higher accuracy'' because they may violate consistency conditions\cite{TAC:linearizable}.  See \cite{Kang_2008_convergence} for details and \cite{cmp:roadmap} for a unified framework. } These indicators of optimality validate that the solution presented in Fig.~\ref{fig:BB-scaledPD} is at least an extremal.

\subsubsection{Illustrating Universality of Proposition $A$}
\begin{figure}[h!]
      \centering
      {\includegraphics[angle=0, width = 0.46\textwidth]{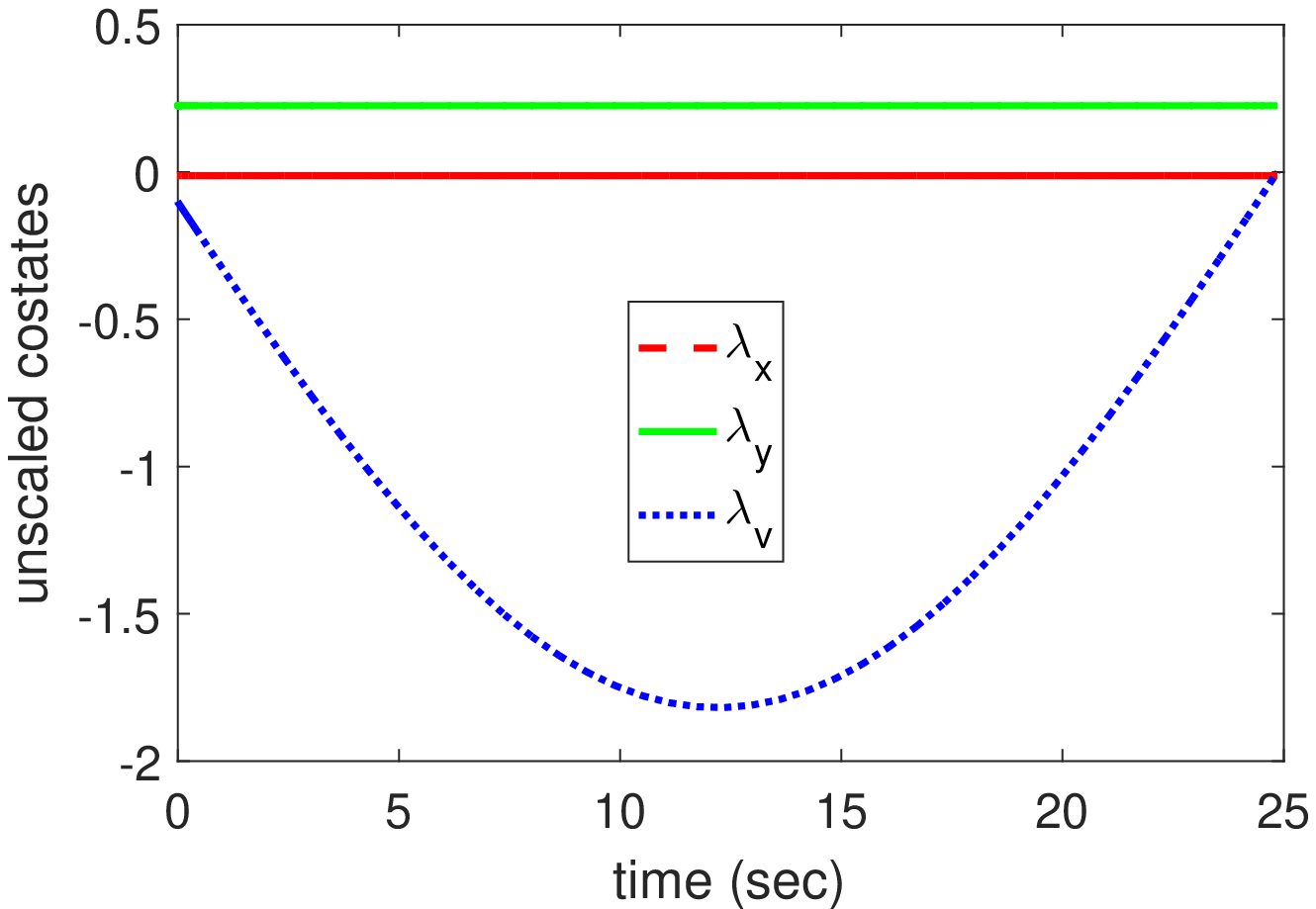}}
            {\includegraphics[angle=0, width = 0.46\textwidth]{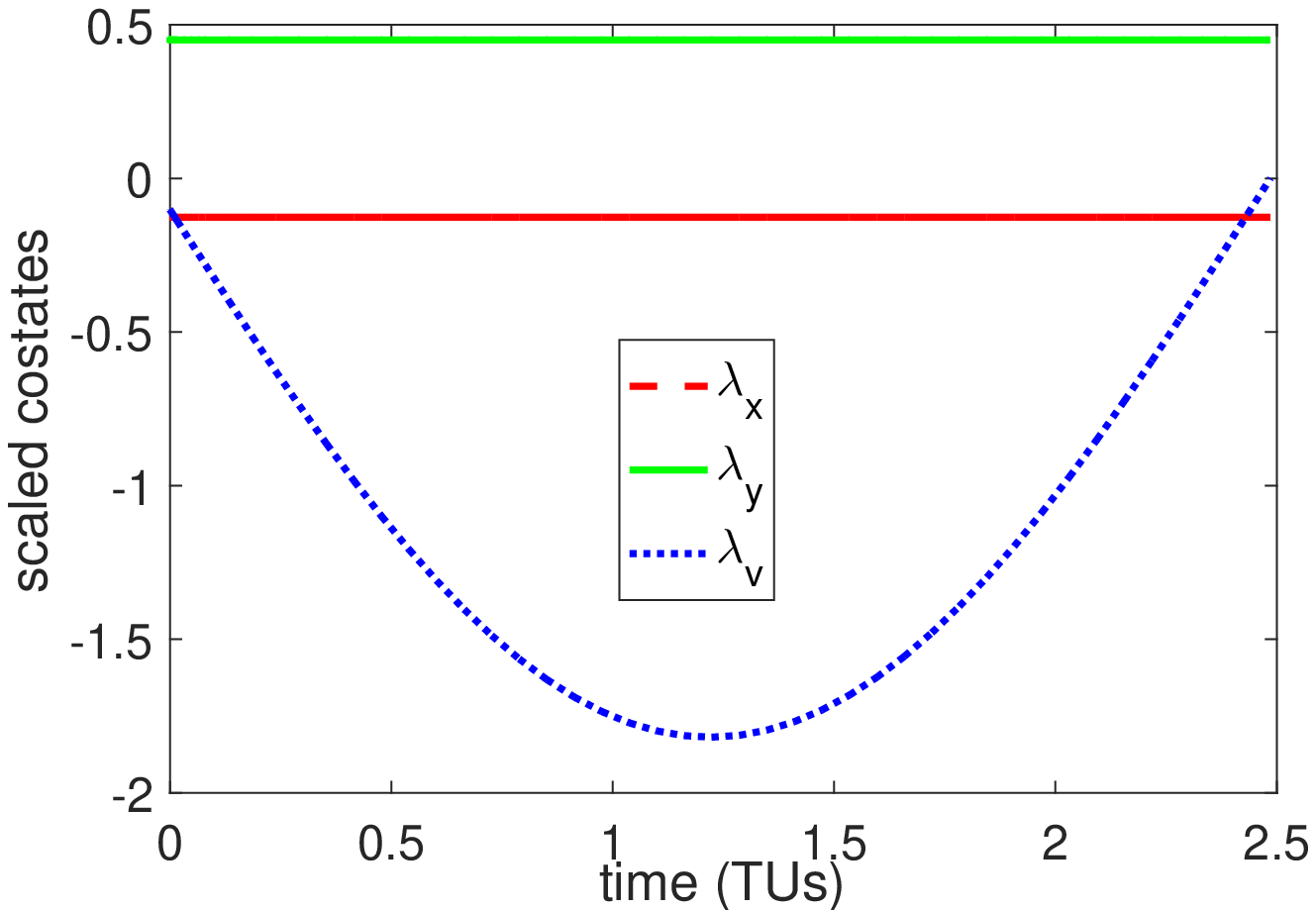}}
      \caption{\textsf{Unscaled and scaled costates obtained by a shooting method for the badly-scaled Brachistochrone problem.}}
    \label{fig:BB-costates-shoot}
   \end{figure}

To illustrate that Proposition $A$ is indeed universal (and not merely specific to DIDO), we construct a shooting algorithm using \textsf{ode45} and \textsf{fsolve} from the MATLAB optimization toolbox.  The objective of \textsf{ode45} is to generate the vector function,
\begin{equation}\label{eq:Bb-Sfun}
\mathbf{S}: (\lambda_{x_0}, \lambda_{y_0}, \lambda_{v_0}, t_f) \mapsto (x_f, y_f, v_f, H_f)
\end{equation}
by integrating the six state-costate equations using the initial conditions (at $t_0 = 0$),
$$ x(t_0) = 0, y(t_0) = 0, v(t_0) = 0, \lambda_x(t_0) = \lambda_{x_0},  \lambda_y(t_0) = \lambda_{y_0}, \lambda_v(t_0) = \lambda_{v_0}   $$
The quantity $H_f$ in \eqref{eq:Bb-Sfun} is the final value of the Hamiltonian evaluated using the results of the integration and \eqref{eq:Bb-Hamil}.  The objective of \textsf{fsolve} is to solve for the zeros of the residual vector function $\br:\real{4} \rightarrow \real{4}$ defined by,
\begin{equation}\label{eq:Bb-gfun}
\br(\lambda_{x_0}, \lambda_{y_0}, \lambda_{v_0}, t_f) := \textbf{S}(\lambda_{x_0}, \lambda_{y_0}, \lambda_{v_0}, t_f) - \left(
                        \begin{array}{c}
                          1000 \\
                          1 \\
                          0 \\
                          -1 \\
                        \end{array}
                      \right)
\end{equation}
Because a shooting algorithm is fundamentally doomed by the curse of sensitivity\cite{ross-book}, we choose the values of the guess to be almost exactly equal to the expected solution,
\begin{equation}\label{eq:Bb-DIDO-init4shoot}
\lambda_{x_0} = -0.013,\quad \lambda_{y_0} = 0.225,\quad \lambda_{v_0} = -0.113,\quad t_f = 24.0
\end{equation}
For the purposes of brevity, we limit our discussions to only the costate trajectories.  Shown in Fig.~\ref{fig:BB-costates-shoot} are both the unscaled and scaled costates trajectories obtained by the shooting algorithm.
%
%
The scaled costates were obtained by using the scaled equations and replacing the numerical value of the 4-vector in \eqref{eq:Bb-gfun} by its scaled counterpart (see \eqref{eq:Bb-xy-final-t}).  It is apparent that the unscaled costates match the descaled costates (see Fig.~\ref{fig:BB-costates}) and the scaled costates match the DIDO result shown in Fig.~\ref{fig:BB-scaledPD} to numerical precision. In other words, we have demonstrated that Proposition \textsl{A} is independent of the numerical algorithm or software.

\subsubsection{Illustrating a Process for Better Balancing}
\begin{figure}[h!]
      \centering
      {\includegraphics[angle=0, width = 0.49\textwidth]{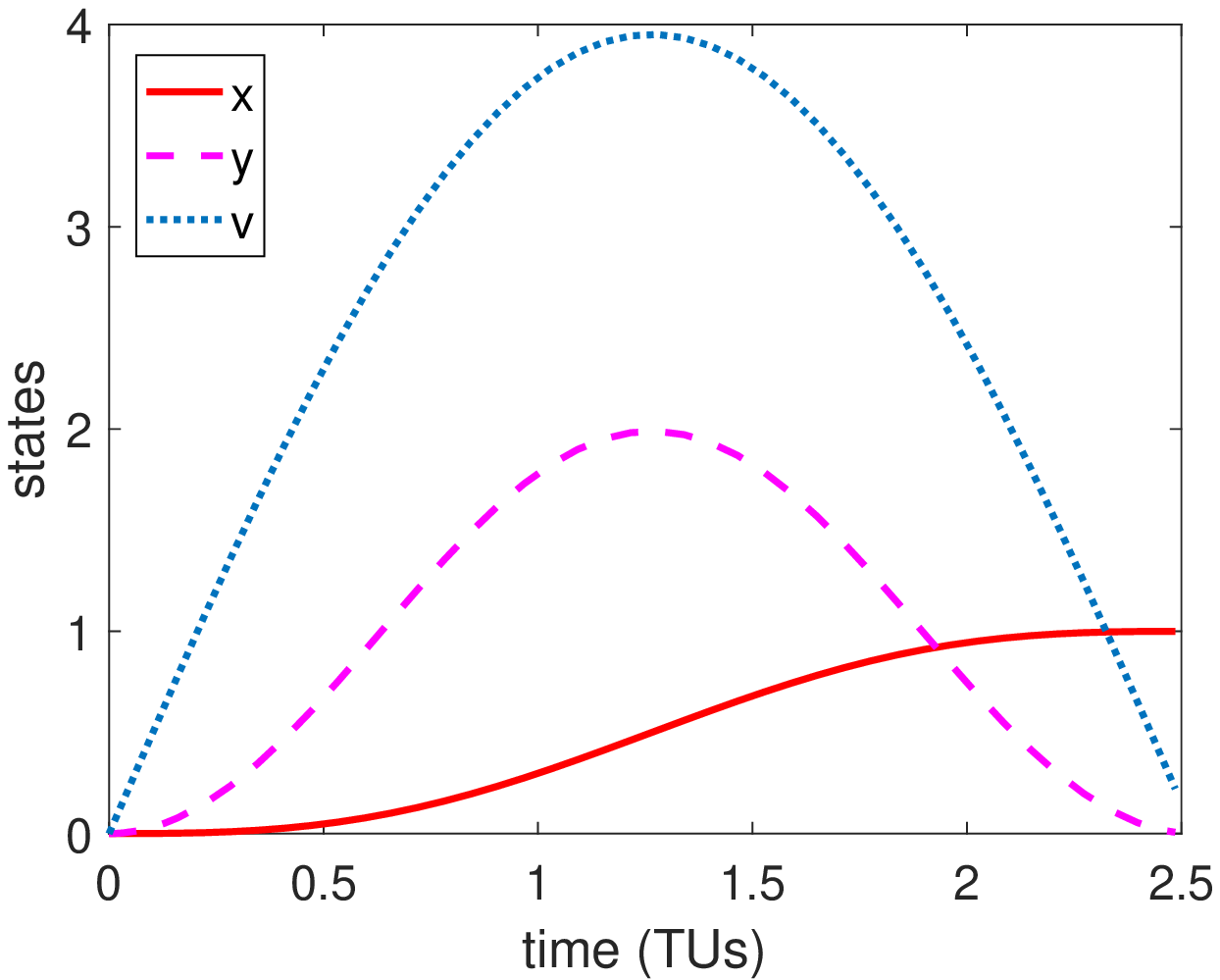}}
            {\includegraphics[angle=0, width = 0.49\textwidth]{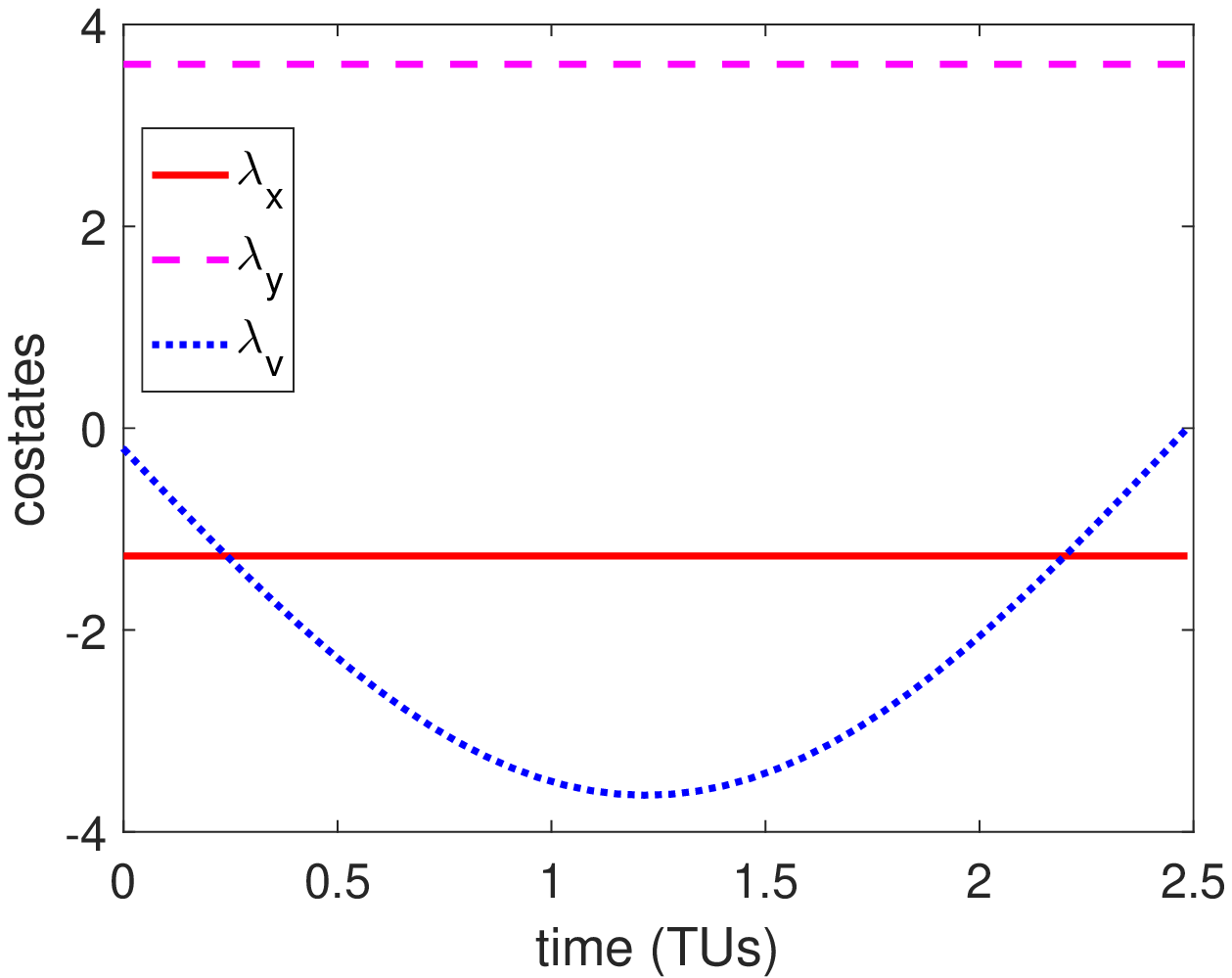}}
      \caption{\textsf{Rescaled solution to the badly-scaled Brachistochrone problem obtained by using the better-balanced set of designer units given by \eqref{eq:BB-designer-units-2}.}}
    \label{fig:BB-RescaledShoot}
   \end{figure}

As a final point of illustration, we now demonstrate that it is possible to achieve a more balanced computational optimal control problem.  First, note from the range of values of the ordinates in Fig.~\ref{fig:BB-scaledPD}, that the computational problem is not perfectly balanced.  This is precisely what happens in solving many flight application problems; that is, it is frequently not possible to choose designer units that achieve well-balanced equations at the first attempt.  Nonetheless, after further analysis of the type illustrated in the preceding paragraphs, it is possible to achieve a better balanced computational problem by merely inspecting the results.  Based on the range of values indicated in Fig.~\ref{fig:BB-scaledPD}, we now rescale the primal problem using the following units:
\begin{equation}\label{eq:BB-designer-units-2}
\begin{aligned}
 x &=  P_x\ \widetilde{x} &&= 1000\ \widetilde{x}  && \Rightarrow \text{1 distance unit along $x$-axis}  && = 1000 \text{ meters} \\
 y & = P_y\ \widetilde{y} &&= 160\ \widetilde{y} && \Rightarrow \text{1 distance unit along $y$-axis}  &&= 160 \text{ meters} \\
 v & = P_v\ \widetilde{v} &&= 20\ \widetilde{v} && \Rightarrow \text{1 speed unit }  &&= 20 \text{ meters/second}
\end{aligned}
\end{equation}
All other choices of units are the same as before; see \eqref{eq:BB-designer-units-1}.  Clearly, $\sqrt{\widetilde{x}^2+\widetilde{y}^2 }$ is meaningless. Note also that the numerical choice of these scaling factors further the disparity between the new set of designer units and the original physical units.  For instance, the gravitational acceleration transforms according to
\begin{equation}\label{eq:g-units-2}
\widetilde{g} = \left(\frac{p_t}{P_v}\right) g  = 4.9\ \widetilde{v}\text{-unit}/\widetilde{t}\text{-unit}
\end{equation}
The primal and dual trajectories for the rescaled problem are shown in Fig.~\ref{fig:BB-RescaledShoot}.
%
%
By inspection, it is clear that the problem is reasonably well-balanced with all variables contained in the range $[-4, 4]$.

\subsubsection{Illustrating the Fallacy of Balancing on the Unit Interval}

\begin{figure}[h!]
      \centering
      {\includegraphics[angle=0, height=0.35\textwidth, width = 0.49\textwidth]{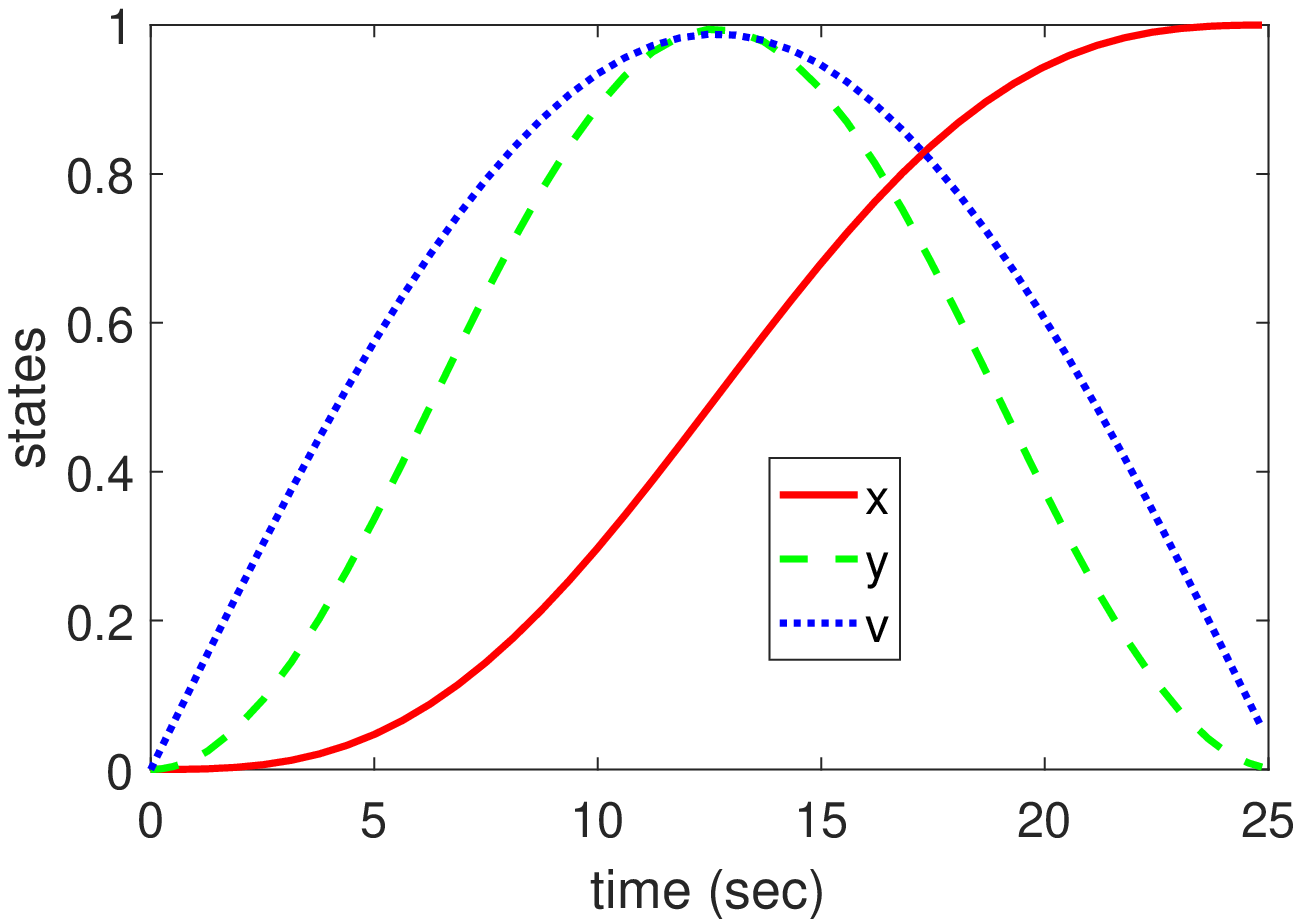}}
            {\includegraphics[angle=0, height=0.35\textwidth, width = 0.49\textwidth]{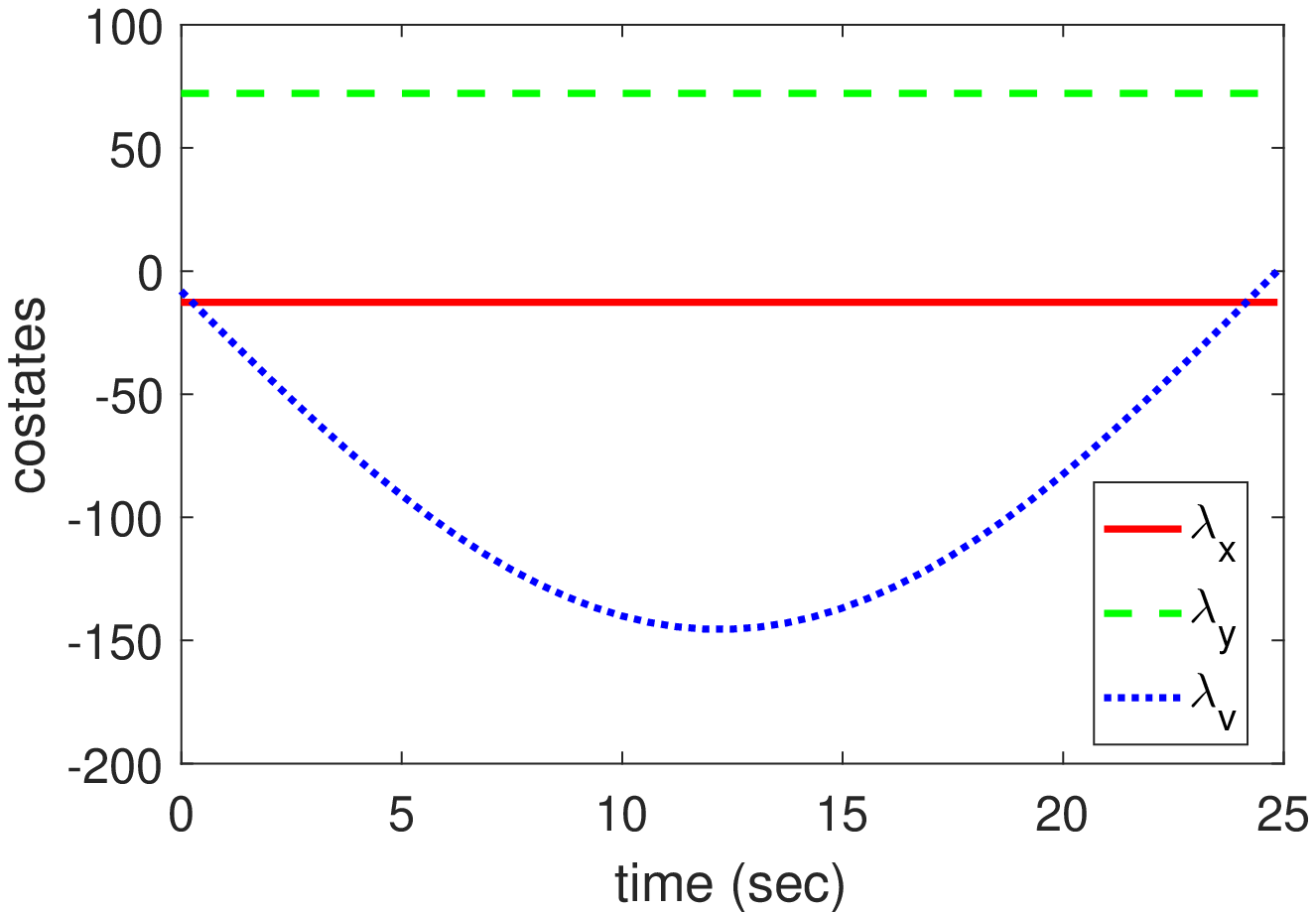}}
      \caption{\textsf{Primal-normalized solution to the badly-scaled Brachistochrone problem obtained by using the set of designer units given by \eqref{eq:BB-designer-units-3}.}}
    \label{fig:BB-RescaledShootUnit}
   \end{figure}

Suppose we scale the primal problem using the following designer units:
\begin{equation}\label{eq:BB-designer-units-3}
\begin{aligned}
 x &=  P_x\ \widetilde{x} &&= 1000\ \widetilde{x}  \\
 y & = P_y\ \widetilde{y} &&= 320\ \widetilde{y}  \\
 v & = P_v\ \widetilde{v} &&= 80\ \widetilde{v}\\
 t & = p_t\ \widetilde{t} && = \widetilde{t}
\end{aligned}
\end{equation}
A quick examination of Fig.~\ref{fig:BB-primalProp} shows that \eqref{eq:BB-designer-units-3}  will force all the state variables to lie on the unit interval $[0, 1]$.  That this is indeed the case is shown in Fig.~\ref{fig:BB-RescaledShootUnit}.
%
%
Also shown in Fig.~\ref{fig:BB-RescaledShootUnit} are the corresponding costates. As expected (by Proposition $A$ and by inspection of  Fig.~\ref{fig:BB-costates}) the dual variables now lie in the range $[-200, 100]$.  In this particular situation, the range $[-200, 100]$ is only about two orders of magnitude away from the desirable interval of $[-1, 1]$; hence, it is no cause for serious alarm. Nonetheless, it is apparent that the imbalance may be greater in other applications.

In certain applications, the state variables may be naturally constrained to lie on the unit interval (e.g., quaternion parametrization). If the costates (e.g., co-quaternions) are imbalanced, it is still possible to balance the state-costate pair by scaling the states to lie on a non-unit interval through a proper selection of $\bP_x$ in Proposition $A$.  In other words, there is no real reason to be constrained on a unit interval.  In stronger terms, \emph{\textbf{ the conventional wisdom of scaling on a unit interval is not necessarily the best approach to balancing optimal control problems}}.

\section{Adjunct Consequences of Scaling and Balancing}

The consequences of scaling and balancing go far beyond faster computation of optimal trajectories. When the pseudospectral optimal control method of the early days\cite{elnagar1,fahroo:LGL1,lncis,fahroo:cheb-jgcd,Williams_2004,Radau-GNC05, AAS:gauss} (c.~1995 - c.~2005) is applied to solve the bad Brachistochrone problem; i.e., without using proper scaling or spectral techniques\cite{knots,spec-alg,arb-grid}, the approach fails to produce the correct solution.
In general, the same is true of Runge-Kutta collocation methods; see \cite{TAC:linearizable} and \cite{paris:PS}. If an algorithm does not generate a feasible solution to a problem with a known solution -- such as the Brachistochrone problem -- then  it is a clear failure of the algorithm.  Frequently, we use optimal control techniques to solve ``hard'' problems where we do not know in advance if a solution exists. In such practical situations it is important to know if the lack of a feasible solution is due to a failure of the algorithm or a genuine non-existence of a solution.  Consequently, over the last decade, proper scaling and balancing have been fundamentally intertwined with the theory and practice of optimal control\cite{conway:survey,ross-book, bhatt:opm,zpm:NASA-report,TEI-JGCD-2011,TRACE-IEEE-Spectrum}.

\subsection{Feasibility via Optimization}

In many practical applications, there is a need to simply generate feasible solutions.
Consequently, scaling and balancing techniques are critical not only for faster optimization but also to answer the more fundamental and difficult theoretical question\cite{vinter} of the existence of a solution.

As a means to illustrate this aspect of the intertwining between theory and practice, consider the optimal propellant maneuver (OPM) that is actively used\cite{bhatt:opm} in current flight operations of the International Space Station (ISS). The OPM saves NASA $90\%$ propellant in momentum management at a savings of approximately \$ 1,000,000 per maneuver\cite{bhatt:opm,SIAMnews,zpm:IEEE}.  This maneuver was designed by Bedrossian after his discovery of the zero propellant maneuver (ZPM) that saves all $100\%$ of propellant\cite{zpm:NASA-report}. Prior to Bedrossian's discovery, it was generally thought that it was impossible to dump all of the accumulated momentum without any propellant consumption\cite{hattis}.  In other words, a feasible solution (for $100\%$ propellant savings) was believed to be nonexistent.  When the ``unscaled'' values of the ISS angular momentum, $h$, the control torque $u$, and angular velocity $\omega$ are used for dynamic optimization, all prior methods investigated by Bedrossian et al consistently failed to generate a feasible solution. These variables vary as\cite{zpm:NASA-report}:
\[
\begin{aligned}
-10^{4} & \le  h \le 10^{4}     && \textit{lb-ft-s}\\
-10^{2}  & \le  u \le 10^{2}          && \textit{lb-ft} \\
-10^{-4} & \le \omega \le 10^{-4}  && s^{-1}
\end{aligned}
\]
Consequently, a failure to find a feasible zero-propellant solution was consistent with the pre-ZPM belief of physics. Were it not for proper scaling and balancing, the ZPM might have gone undiscovered. More specifically, when the variables are transformed according to
\begin{equation}\label{eq:ZPM-transformation-1}
\begin{aligned}
h &= 1000\ \widetilde{h}     && \textit{lb-ft-s}\\
u  & = 10\  \widetilde{u}          && \textit{lb-ft} \\
t & = 1000\ \widetilde{t}   && s
\end{aligned}
\end{equation}
not only were Bedrossian et al\cite{zpm:NASA-report} able to find a feasible solution but also several different solutions! In other words, a special choice of designer units ``converted'' a hard problem to an easy one. The rest is history\cite{SIAMnews}.

The lessons learned from such successes and similar ones that followed\cite{bhatt:opm,TRACE-IEEE-Spectrum, Kepler-micro-slew,Karp-JWST,CDC-Workshop} are codified in the scaling and balancing techniques presented in the preceding sections. More specifically, the last decade has witnessed the use theoretical optimization principles to determine the practical feasibility of innovative concepts as opposed to the more conventional use of algorithms to optimize a feasible design.  In other words, optimal control theory has been used a tool to innovate and not merely to optimize.

\subsection{Nonuniqueness of the Costates}

Note that Proposition \textsl{A} only asserts the existence of linearly descaled covectors; it does not imply uniqueness.  Using the same arguments as that of Proposition \textsl{A}, it is relatively straightforward to show that the costates in an optimal control problem are not necessarily unique.  This is achieved by replacing the linear equations \eqref{eq:result-dual-lam}-\eqref{eq:result-dual-nu} with its affine counterparts:
\begin{subequations}\label{eq:nonunique}
\begin{align}
\blam^*(\cdot) := & \bP_\lambda\ \tlam^*(\cdot) + \bq_\lambda(\cdot)\\
\bmu^*(\cdot) := & \bP_\mu\ \tmu^*(\cdot) + \bq_\mu(\cdot) \\
\bnu^* := & \bP_\nu\ \tnu^* + \bq_\nu
\end{align}
\end{subequations}
where, $\bq_\lambda(\cdot): t \mapsto \real{N_x}$, $\bq_\mu(\cdot) : t \mapsto \real{N_h}$ and $\bq_\nu \in \real{N_e}$.  Obviously, $\bq_\lambda(\cdot) \equiv \bzero  $, $\bq_\mu(\cdot) \equiv \bzero$ and $\bq_\nu = \bzero$ recovers Proposition \textsl{A}; however, by substituting \eqref{eq:nonunique} in Problem $B^\lambda$, it is straightforward to show that it is not necessary for the $\bq$-multipliers to be trivial.  That Lagrange multipliers are not unique is well-known in nonlinear programming\cite{kyparisis,wachsmuth}; hence, it seems apparent that this must also be true in optimal control programming.  However, unlike the static case, the possibility of nonunique costates in optimal control is limited by the Lipschitz-continuity of $\partial_x\bff$. Despite this limitation, the conditions for nonuniqueness are relatively mild; see Ref.~\cite{ross-book}, Sec.~4.9, for a complete worked-out example pertaining to the optimal steering of a rigid body.

\subsection{Fallacy of Discrete Scaling}
In \eqref{eq:transform-var} and \eqref{eq:transform-fun} we deliberately used affine scaling with constant coefficients.  Suppose we choose time-varying scaling coefficients; then, the state variable transformation can be written as,
\begin{equation}\label{eq:transform-var-timeVarying}
\bx(t) =  \bP_x(t)\ \bxt(t) + \bq_x(t)
\end{equation}
Differentiating \eqref{eq:transform-var-timeVarying} and substituting the unscaled dynamics in the resulting equation generates,
\begin{equation}\label{eq:add-dynamics}
\begin{aligned}
\frac{d\bxt}{d\widetilde{t}}  &=  p_t\bP^{-1}_x(t) \left(\dot\bx(t) -  \dot\bP_x(t)\ \bxt(t) - \dot\bq_x(t) \right) \\
&=  p_t\bP^{-1}_x(t) \bff(\bx(t), \bu(t), t) - \underbrace{p_t\bP^{-1}_x(t)\left( \dot\bP_x(t)\ \bxt(t) + \dot\bq_x(t) \right)}_{\text{additional dynamics}}
\end{aligned}
\end{equation}
That is, the transformed state dynamics contain additional dynamics. Consequently, it generates the following questions:
\begin{enumerate}
\item How do we choose $\dot\bP_x(t)$ and $\dot\bq_x(t)$?  That is, what is the rationale for choosing these functions?
\item Because our objective is to generate a theory for generic problems (e.g., generic dynamics) how can we choose universal functions $\dot\bP_x(t)$ and $\dot\bq_x(t)$?
\item Even if we were to severely limit time-varying scaling to a specific dynamical system, how do we choose $\dot\bP_x(t)$ and $\dot\bq_x(t)$ whose properties remain valid for all feasible control functions $\bu(\cdot)$?
\item How do we ensure that the additional dynamics indicated in \eqref{eq:add-dynamics} does not create new numerical problems over the space of all differentiable functions $\bff$?
\end{enumerate}
From these basic considerations, it is clear that time-varying affine scaling generates more questions than answers.

Interestingly, time-varying scales are implicit in many software packages and algorithms. To appreciate this point, consider the discretization of a one-dimensional state trajectory, $t \mapsto x \in \Real$.  For $k = 0, \ldots, N$, the discretized variables, $x_k$ represent the samples of the state trajectory; hence, we can write,
\begin{equation}\label{eq:disc2cont}
\left[
  \begin{array}{c}
    x_0 \\
    x_1 \\
    \vdots \\
    x_N \\
  \end{array}
\right] = \left[
  \begin{array}{c}
    x(t_0) \\
    x(t_1) \\
    \vdots \\
    x(t_N) \\
  \end{array}
\right]
\end{equation}
If the discretized variables are scaled by, say, a diagonal matrix with entries $P_0, \ldots, P_N$, then \eqref{eq:disc2cont} transforms according to
\begin{equation}\label{eq:disc-scaling-1}
\left[
  \begin{array}{c}
    x_0 \\
    x_1 \\
    \vdots \\
    x_N \\
  \end{array}
\right] := \left[
  \begin{array}{c}
   P_0\ \widetilde{x}_0 \\
   P_1\ \widetilde{x}_1 \\
    \vdots \\
   P_N\ \widetilde{x}_N \\
  \end{array}
\right] = \left[
  \begin{array}{c}
    x(t_0) \\
    x(t_1) \\
    \vdots \\
    x(t_N) \\
  \end{array}
\right]
\end{equation}
Let, $P(t)$ be any function such that $P(t_k) = P_k, k = 0, \ldots, N$; then, we can write \eqref{eq:disc-scaling-1} as
\begin{equation}\label{eq:disc-scaling}
\left[
  \begin{array}{c}
    x_0 \\
    x_1 \\
    \vdots \\
    x_N \\
  \end{array}
\right] := \left[
  \begin{array}{c}
   P_0\ \widetilde{x}_0 \\
   P_1\ \widetilde{x}_1 \\
    \vdots \\
   P_N\ \widetilde{x}_N \\
  \end{array}
\right] := \left[
  \begin{array}{c}
    P(t_0)\ \widetilde{x}(t_0) \\
    P(t_1)\ \widetilde{x}(t_1) \\
    \vdots \\
    P(t_N)\ \widetilde{x}(t_N) \\
  \end{array}
\right] = \left[
  \begin{array}{c}
    x(t_0) \\
    x(t_1) \\
    \vdots \\
    x(t_N) \\
  \end{array}
\right]
\end{equation}
Hence, it follows that scaling the discretized variables is equivalent to discretizing the continuous-time trajectory according to,
\begin{equation}\label{eq:disc2cont-insight}
x(t) = P(t) \widetilde{x}(t)
\end{equation}
Consequently, any algorithm or software package that scales the variables at the discrete-level is implicitly using time-varying scales at the optimal-control level. Given this fact, it is critical that the additional dynamics noted in \eqref{eq:add-dynamics} be automatically incorporated at the discrete level in the algorithm or software package. To understand how this can be done, let $\triangle_k$ be any discrete derivative.
Because the additional dynamics in \eqref{eq:add-dynamics} is a consequence of the product rule of continuous calculus, the ``equivalent'' discrete product rule,
\begin{equation}\label{eq:disc-product-rule}
\triangle_k x_k := P_k\left(\triangle_k\, \widetilde{x}_k\right) + \left(\triangle_k P_k\right) \widetilde{x}_k, \quad k = 0, \ldots, N
\end{equation}
must be naturally incorporated in the algorithm or software package.   If $\triangle_k$ is a forward difference operator, then,
$$\triangle_k(P\cdot \widetilde{x})_k := P_{k+1}\widetilde{x}_{k+1} - P_k \widetilde{x}_k $$
It is quite straightforward to show that $\triangle_k(P\cdot \widetilde{x})_k \ne \triangle_k\,x_k$, and that,
\begin{equation}\label{eq:error}
\triangle_k\,x_k - \triangle_k(P\cdot \widetilde{x})_k = (\triangle_k P_k) (\triangle_k \widetilde{x}_k)
\end{equation}
The right-hand-side of \eqref{eq:error} is not necessarily a second-order effect unless the scaling algorithm renders $\abs{(\triangle_k P_k) (\triangle_k \widetilde{x}_k)}$ small.
Recall that $P_k$ is a scaling factor at the discretized level and not necessarily connected to some continuous function $P(t)$ with a small Lipschitz constant; see \eqref{eq:disc-scaling-1}.
\emph{\textbf{Consequently, if $\abs{(\triangle_k P_k) (\triangle_k \widetilde{x}_k)}$ is not small, then the algorithm -- any algorithm -- is attempting to solve for the wrong dynamics!}}  The implications of this insight are far reaching:
\begin{enumerate}
\item If the algorithmic iterations do not converge and/or are expensive (i.e., take long computational time), it is quite possible the original optimal control problem might have been easy but rendered hard because of scaling at the discretized level!
\item The situation might be made even worse with more sophisticated scaling like adaptive scaling. In such schemes, the scale factors $P_0, P_1, \ldots, P_N$ (see \eqref{eq:disc-scaling}) are not ``constants'' over the course of the iteration but change ``adaptively'' based on the current iterate.  In following the same process that led to \eqref{eq:disc2cont-insight}, adaptive scaling implies that in continuous-time the state variable is scaled according to some feedback process,
\begin{equation}
x(t) = P(t, x) \widetilde{x}(t)
\end{equation}
over the course of the iterations. Ignoring the possible stability issues resulting from this feedback process, it is clear that this may be worse than time-varying scaling because we now have additional dynamics associated with $\dot P$,
\begin{equation}
\frac{dP}{dt} = \frac{\partial P}{\partial t} + \frac{\partial P}{\partial x} \frac{dx}{dt}
\end{equation}
These continuous-time dynamics are not necessarily incorporated in adaptive scaling because the chain rule (of continuous calculus) must also be incorporated at the discrete level (in addition to the product rule).
\item From the preceding point, it is also clear that nonlinear scaling also has the same drawbacks of adaptive scaling (at the discrete level).
\item In an optimistic scenario, is quite possible that $\abs{(\triangle_k P_k) (\triangle_k \widetilde{x}_k)}$ is small either by an implicit/explicit result of a scaling algorithm or by accident. Even under this fortuitous case, the error in the satisfaction of the dynamical equations is higher than the computational tolerances enforced unless of course $\abs{(\triangle_k P_k) (\triangle_k \widetilde{x}_k)} = 0$ for $k = 0, \ldots, N$.  \emph{\textbf{This is yet another reason why an independent verification of feasibility as highlighted in Sec.~VI.C is crucial for validating numerical accuracy}}.
\item In the best-case scenario, $\abs{(\triangle_k P_k) (\triangle_k \widetilde{x}_k)} = 0$ for $k = 0, \ldots, N$.  In this case, an algorithm is solving for the correct but transformed dynamics given by \eqref{eq:add-dynamics}.  In the absence of new analysis, there is no apparent reason why the transformed dynamics is universally better for optimization than the original dynamics.
\end{enumerate}
The preceding analysis explains why autoscaling done at the discrete level without explicit consideration of the dynamics of the optimal control problem may be harmful to the accuracy and convergence of the algorithm.
A simple remedy for this problem is to perform scaling and balancing at the optimal-control level, and turn off any autoscaling options of software packages that are based on scaling the discretized Jacobian that does not incorporate the additional dynamics presented in \eqref{eq:add-dynamics}.

\section{Conclusions}

Some of the ideas and parts of the process presented in this paper have been used for well over a decade -- albeit in an application-specific manner -- to generate successful flight implementations. In this paper, we have generalized previous concepts and provided a clear mathematical framework leading to a more unified procedure.  More specifically, we have shown that the concept of designer units is fundamentally liberating. If necessary, it even allows one to choose radically different units of measurement along $x$- and $y$-directions. As a consequence of this liberation, the physical concept of a vector as quantity with magnitude and direction must be abandoned.  A vector is simply a stack of scalar variables in any units.  A covector is a measurement conversion device that connects the disparate units of a vector to some common unit.  In an optimal control problem, the Lagrange multipliers are covectors, and the common unit of measurement is the cost unit or the cost-unit per time-unit. \emph{\textbf{The numerical values of a covector have a seesaw effect on the values of the associated vector. The seesaw effect can be used to balance the primal and dual variables for computational efficiency}}.

It is not necessary to scale and balance the variables on a unit interval; in fact, it may not even be feasible.  Even in situations where the state variables are naturally constrained to the unit interval, it is possible to scale them beyond their physical bounds to achieve better balancing.

Scaling an optimal control problem at the nonlinear-programming level is likely to induce unwanted dynamics that may render an easy problem hard. Therefore, great caution must be exercised in using software packages that simply patch discretization methods to nonlinear programming solvers. In contrast, when scaling and balancing is done at the optimal-control level, it can be used quite powerfully to innovate by solving hard problems easily.


\section{Acknowledgments}

We thank Scott Josselyn (former graduate student) who showed us (nearly two decades ago!) that optimal control problems may be scaled ``inconsistently'' to achieve high computational efficiency. We also thank Steve Paris (formerly of the Boeing company) for providing the initial spark of using Lagrange multipliers as ``balancing'' devices.  Last but not least, we thank the anonymous reviewers for providing many constructive comments which greatly improved the quality of this paper.

\end{document}